\documentclass[11pt]{amsart}

\usepackage{amssymb,amsthm,amsmath}
\usepackage[numbers,sort&compress]{natbib}
\usepackage[margin=1in]{geometry}
\usepackage{color}
\usepackage{graphicx}
\usepackage{tikz}
\usepackage[colorlinks,
linkcolor=red,
anchorcolor=blue,
citecolor=green
]{hyperref}

\usepackage{enumerate}

\def\e{\varepsilon}

\let\newpf\proof \let\proof\relax 
\newenvironment{pf}{\newpf[\proofname]}{\qed\endtrivlist}

\newcommand{\ba}{\overline{A}}

\def\be{\begin{equation}}
\def\ee{\end{equation}}

\def\ba{{\begin{align}}}
\def\ea{{\end{align}}}

\def\bm{\begin{matrix}}
\def\em{\end{matrix}}

\def\0{{\mathbf 0}}

\newtheorem{Theorem}{Theorem}[section]
\newtheorem{Lemma}{Lemma}[section]
\newtheorem{Proposition}{Proposition}[section]
\newtheorem{Corollary}{Corollary}[section]
\theoremstyle{definition}
\newtheorem{Remark}{Remark}[section]
\newtheorem{Example}{Example}[section]
\newtheorem{Definition}{Definition}[section]

\numberwithin{equation}{section}

\theoremstyle{definition}

\renewcommand{\mod}{\operatorname{mod}}

\newcommand{\C}{{\mathbb C}}

\newcommand{\N}{{\mathbb N}}
\newcommand{\Q}{{\mathbb Q}}
\newcommand{\R}{{\mathbb R}}
\newcommand{\T}{{\mathbb T}}

\newcommand{\Z}{{\mathbb Z}}

\catcode`\@=12

\def\Empty{}
\newcommand\oplabel[1]{
  \def\OpArg{#1} \ifx \OpArg\Empty {} \else
    \label{#1}
  \fi}

\newcommand{\comm}[1]{}
\newcommand{\comment}[1]{}

\begin{document}

\title[]{Transition space for the continuity of the Lyapunov exponent of quasiperiodic Schr\"odinger cocycles}

\author{Lingrui Ge}
\address{Department of Mathematics, University of California Irvine, CA, 92697-3875, USA}
\email{lingruig@uci.edu}

\author{Yiqian Wang}
\address{
Department of Mathematics, Nanjing University, Nanjing 210093, China}
\email{yiqianw@nju.edu.cn}

\author {Jiangong You}
\address{
Chern Institute of Mathematics and LPMC, Nankai University, Tianjin 300071, China} \email{jyou@nankai.edu.cn}

\author{Xin Zhao}
\address{
Department of Mathematics, Nanjing University, Nanjing 210093, China
and
Department of Mathematics, University of California Irvine, CA, 92697-3875, USA}
 \email{njuzhaox@126.com}

\begin{abstract}
We construct discontinuous point of the Lyapunov exponent of quasiperiodic Schr\"odinger cocycles in the Gevrey space $G^{s}$ with $s>2$. In contrast, the Lyapunov exponent has been proved to be continuous  in the Gevrey space $G^{s}$ with $s<2$ \cite{klein,cgyz}. This shows that  $G^2$ is the transition space for the continuity of the Lyapunov exponent.
\end{abstract}

\maketitle
\section{Introduction}
Let $X$ be a $C^r$ compact manifold, $A(x)$ be a $SL(2,\R)$-valued function on $X$ and $(X,T,\mu)$ be ergodic with $\mu$ a normalized $T$-invariant measure. Dynamical systems on $X \times \R^2$  given by
$$
(x,w) \rightarrow (T(x),A(x)w)
$$
are called a $SL(2,\R)$-cocycle and  denoted by $(T,A)$. In particular, if $X=\mathbb{T}^n=\R^n/2\pi\Z^n$ and $T=T_\alpha: x\rightarrow x+2\pi\alpha$ with $\alpha$ independent over $\Q$, we call $(T_\alpha,A)$ a quasiperiodic $SL(2,\R)$-cocycle, which is simply denoted by $(\alpha,A)$. If moreover,
$$
A(x)=S_{E,v}(x)=\begin{pmatrix}
E-v(x)& -1\\
1& 0
\end{pmatrix}
$$
with $v(x)$ a $2\pi$-periodic function in each variable, we call $(\alpha,S_{E,v})$ a quasiperiodic Schr\"odinger cocycle.

The $n$-th iteration of the cocycle $(T,A)$ is denoted by $(T,A)^n=(T^n,A_n)$ where
\begin{align*}
\begin{split}
A_n(x)=\left\{
\begin{array}{ll}
A(T^{n-1}x)\cdots A(x), & n\geq 1\\
I_2, & n= 0\\
A_{-n}(T^nx)^{-1}, & n\leq -1
\end{array}
\right.
\end{split}.
\end{align*}
The (maximum) Lyapunov exponent $L(A)$ of the cocycle is defined as
$$
L(A)=\lim\limits_{n\rightarrow\infty}\frac{1}{n}\int_{X}\ln\|A_n(x)\|d\mu
=\inf_n\frac{1}{n}\int_{X}\ln\|A_n(x)\|d\mu\geq 0.
$$
The limit
exists and is equal to the infimum since $\left\{\int_{X}\ln\|A_n(x)\|
d\mu\right\}_{n\geq1}$ is a subadditive sequence. Moreover, by Kingman's subadditive ergodic theorem, we also
have
$$L(A)=\lim\limits_{n\rightarrow\infty}\frac{1}{n}\ln\|A_n(x)\|$$
for $\mu$-almost every $x\in X$.

Regularity of the Lyapunov exponent (LE) is one of the central subjects in
smooth dynamical systems, which depends subtly on the base dynamics $T$ and the smoothness of the matrix $A$. In the present paper, we are mainly interested in how the regularity of $A$ affects the continuity of LE for quasiperiodic $SL(2,\R)$/Schr\"odinger cocycles. Our motivation comes from the pioneering (opposite) results on the continuity of the Lyapunov exponent in $C^\omega$ and $C^\infty$ spaces.
\begin{itemize}
\item For any quasiperiodic $SL(2,\R)$-cocycle, the Lyapunov exponent is always continuous with respect to  $A$ in  $C^\omega$ topology \cite{bj,ajs}.
\item  the Lyapunov exponent is not always continuous with respect to  $A$ in  $C^\infty$ topology \cite{wangyou,wy2}.
\end{itemize}
It is well-known that  the Gevrey spaces $G^s, 1<s<\infty$ are between the $C^\infty$ and analytic spaces. Roughly speaking,  a $2\pi$-periodic cocycle map $A\in G^s$ means that  $\widehat{A}(k)$, the Fourier coefficients of $A$, decay sub-exponentially like $\mathcal{O}(e^{-k^{1/s}h})$, while the Fourier coefficients of an analytic function decay exponentially and the Fourier coefficients of a smooth function decay faster than any polynomials. In this paper, we are interested in finding the optimal Gevrey space  to ensure the continuity of the Lyapunov exponent. More concretely,  we prove that
the  Lyapunov exponent of quasiperiodic Schr\"odinger cocycles is discontinuous in the Gevrey space $G^{s}$ with $s>2$. In contrast, the Lyapunov exponent is continuous  in $G^{s}$ with $s<2$ \cite{klein,cgyz}. This shows that  $G^2$ is the transition space for the continuity of the Lyapunov exponent.

It is known that a powerful tool to prove the continuity of LE is the large deviation theorem (LDT) and avalanch principle (AP). Our results in some sense show that LDT breaks down for general $G^s(\mathbb{S},SL(2,\R))$-cocycles with $s>2$. One can also compare our result with the result in \cite{klein} where Klein showed LE is continuous with respect to the energies E if the potential is in an open and dense subspace of $G^s(\mathbb{S}^1)$ with $s>2$ where certain non-degeneracy condition is satisfied. Our result shows that such non-degeneracy condition is necessary as  the LE would be discontinuous if the potential is sufficiently ``degenerate".

We finally remark that transition phenomenon seems to be a common phenomenon in quasiperiodic dynamical systems and always attracts people's attention. For example, it was shown in \cite{cw} that
any Lagrangian torus with a given unique rotation vector  of an integrable Hamiltonian can be destructed by an
arbitrarily $C^{2d-\delta}$-small perturbation. In contrast, it was shown that
KAM torus with constant type frequency persists for all $C^{2d+\delta}$-small perturbations
\cite{poschel}. Thus $C^{2d}$ is the transition space for the persistence of KAM torus.
\subsection{Transition phenomena for quasiperiodic Schr\"odinger operators}

 The discrete one dimensional quasiperiodic Schr\"odinger operators on $\ell^2(\Z)$ are given by
\begin{equation}\label{schrodinger}
(H_{\lambda v,x,\alpha}u)_n=u_{n+1}+u_{n-1}+\lambda v(x+n\alpha)u_n,\ \ n\in\Z,
\end{equation}
where $\alpha\in\R^d$ is the frequency, $x\in\mathbb{S}^d$ is the phase, $\lambda\in\R$ is the coupling constant and $v\in C^r(\mathbb{S}^d,\R)$ ($r=0,1,\cdots,\infty,\omega$) is called the potential. The spectral properties of operator \eqref{schrodinger} is closely related to the Schr\"odinger cocycle $(\alpha, S_
{E,\lambda v})\in \mathbb{S}^d\times C^r(\mathbb{S}^d, SL(2,\R))$.
Quasiperiodic Schr\"odinger operators naturally arise in solid-state physics, describing the influence of an external magnetic field on the electrons of a crystal.

Different from random Schr\"odinger operators, an important feature of one dimensional quasiperiodic operators is that the family $\{H_{\lambda v,x,\alpha}\}_{\lambda\in\R}$ undergoes a so called metal-insulator transition when $|\lambda|$ changes from small to large. Indeed, besides the metal-insulator transition, various spectral transition phenomena take place for quasiperiodic operators. Here we give some perfect examples.
\begin{Example}[Metal-insulator transition]Assume $\alpha$ is Diophantine \footnote{ $\alpha \in\R$ is  called {\it Diophantine}, denoted by $\alpha \in {\rm DC}(\kappa,\tau)$, if there exist $\kappa>0$ and $\tau>1$ such that
\begin{equation}\label{dio}
{\rm DC}(\kappa,\tau):=\left\{\alpha \in\R:  \inf_{j \in \Z}\left|n\alpha- j \right|
> \frac{\kappa}{|n|^{\tau}},\quad \forall \  n\in\Z\backslash\{0\} \right\}.
\end{equation}
} and $v(x)=2\cos2\pi x$, the following results were given by Jitomirskaya \cite{J} in 1999,
\begin{itemize}
\item $|\lambda|>1$, $H_{\lambda v,x,\alpha}$ has Anderson localization for a.e. $x$,
\item $|\lambda|=1$, $H_{\lambda v,x,\alpha}$ has purely singular continuous for a.e. $x$,
\item $|\lambda|<1$, $H_{\lambda v,x,\alpha}$ has purely absolutely continuous spectrum for a.e. $x$.
\end{itemize}
\end{Example}

\begin{Example}[Sharp spectral transition in frequency] We denote by
\begin{equation*}
\beta(\alpha)=\limsup\limits_{k\rightarrow\infty}-\frac{\ln\|k\alpha\|_{\R/\Z}}{|k|},
\end{equation*}
where $\|x\|_{\R/\Z}=\text{dist}(x,\Z).$ Let $v(x)=2\cos2\pi x$, i.e., the famous almost Mathieu operators,
\begin{itemize}
\item $|\lambda|>e^{\beta(\alpha)}$, $H_{\lambda v,x,\alpha}$ has Anderson localization for a.e. $x$ \cite{ayz,JLiu} \footnote{\cite{ayz} proved the measure version and \cite{JLiu} proved the arithmetic version, actually \cite{JLiu} proved Anderson localization holds for Diophantine phases.},
\item $1\leq |\lambda|<e^{\beta(\alpha)}$, $H_{\lambda v,x,\alpha}$ has purely singular continuous for all $x$ \cite{ayz,J1} \footnote{\cite{ayz} proved the case $\beta>0$, \cite{J1} proved the case $\alpha$ irrational and $|\lambda|=1$.},
\item $|\lambda|<1$, $H_{\lambda v,x,\alpha}$ has purely absolutely continuous spectrum for all $x$ \cite{Avila1}.
\end{itemize}
\end{Example}

\begin{Example}[Sharp spectral transition in phase] We denote by \begin{equation*}
\delta(\alpha,x)=\limsup\limits_{k\rightarrow\infty}-\frac{\ln\|2x+k\alpha\|_{\R/\Z}}{|k|}.
\end{equation*}
Let $v(x)=2\cos2\pi x$, i.e., the famous almost Mathieu operators,
\begin{itemize}
\item $|\lambda|>e^{\delta(\alpha,x)}$, $H_{\lambda v,x,\alpha}$ has Anderson localization for Diophantine $\alpha$ \cite{JLiu1},
\item $1\leq |\lambda|<e^{\delta(\alpha,x)}$, $H_{\lambda v,x,\alpha}$ has purely singular continuous for all $\alpha$ \cite{JLiu1,J1}.
\item $|\lambda|<1$, $H_{\lambda v,x,\alpha}$ has purely absolutely continuous spectrum for all $\alpha$ \cite{Avila1}.
\end{itemize}
\end{Example}

 Although the transition phenomenon is common for quasiperiodic Schr\"odinger operators, however, the exact transition points are usually difficult to obtain as it depends sensitively on the arithmetic properties of the frequency and phase.
 This paper will give   explicit transition  space for the continuity of the Lyapunov exponent. For our purpose, we introduce the following space of Gevrey functions and its topology.

  For any smooth function $f$ defined on $\mathbb{S}^1$, let
$$
|f|_{s,K}:=\frac{4\pi^2}{3}\sup\limits_{k}\frac{(1+|k|)^2}{K^k(k!)^s}|\partial^k f|_{C^0(\mathbb{S}^1)},
$$
$G^{s,K}(\mathbb{S}^1)=\{f\in C^\infty(\mathbb{S}^1,\R)||f|_{s,K}<\infty\}$ and $G^s(\mathbb{S}^1)=\cup_{K>0}G^{s,K}(\mathbb{S}^1)$. Note that $G^{s,K}(\mathbb{S}^1)$ is a Banach space. Obviously, $G^1(\mathbb{S}^1)$  is the space of analytic functions and for any $s\ge 1$, $G^s(\mathbb{S}^1)$ is a subspace of the space of smooth functions. We equip $G^s(\mathbb{S}^1)$ with the usual inductive limit topology. That is, $f_n$ converges to $f$ in $G^{s}(\mathbb{S}^1)$-topology if and only if  $|f_n-f|_{s,K}\rightarrow 0$ as $n\rightarrow \infty$ for some $K>0$.

We say $\alpha\in\R\backslash \Q$ is of bounded type if there exists $M>0$, such that the continued fraction expansion of $\alpha$, denoted by $p_n/q_n$ satisfying
$$
q_{n+1}\leq Mq_n,\ \ \forall n\in\N.
$$

\begin{Theorem}\label{main}
Assume $\alpha$ is of bounded type, for quasiperiodic Schr\"odinger cocycle $(\alpha, S_{E,v})$, we have
\begin{enumerate}[$(1)$]
\item For any $v\in G^s(\mathbb{S}^1)$ with $s<2$, the Lyapunov exponent is continuous with respect to $v$ in $G^s$-topology.
\item There exists $v_0\in G^s(\mathbb{S}^1)$ with $s>2$, such that the Lyapunov exponent is discontinuous at $v_0$ in $G^s$-topology.
\end{enumerate}
\end{Theorem}
\begin{Remark}
Part (1) of Theorem \ref{main} was recently proved in \cite{cgyz} (Theorem 6.3 of \cite{cgyz}), we list here for completeness. The main result of the present paper is part (2).
\end{Remark}
\begin{Remark}
The bounded type $\alpha$ is dense in $\R$.
\end{Remark}
Part (2) of Theorem \ref{main} can be obtained in the same way as in \cite{wangyou} (See page 2367, proof of Theorem 2 in \cite{wangyou}) from the following examples in $SL(2,\R)$-cocycles.
\begin{Theorem}\label{main1}
Consider quasiperiodic $SL(2,\R)$-cocycles over $\mathbb{S}^1$ with $\alpha$ being a fixed
irrational number of bounded-type. For any $s>2$, there exists a cocycle $D_s\in G^s(\mathbb{S}^1, SL(2,\R))$ such that the Lyapunov exponent is discontinuous at $D_s$ in $G^s$-topology.
\end{Theorem}
\subsection{A Brief review on the continuity of the Lyapunov exponent}
As we mentioned above, both the base dynamics $T$ and the smoothness of the matrix $A$ affect the regularity of the Lyapunov exponent. This has been the object of considerable recent interests, see Viana \cite{viana}, Wilkinson \cite{wilk} and the references therein.

If the base dynamics has some hyperbolicity, then the Lyapunov exponent is continuous. For example, Furstenberg-Kifer \cite{FK}
and Hennion \cite{H} proved the continuity of the largest LE of i.i.d random matrices under a
condition of almost irreducibility. Bocker and Viana \cite{BV} proved continuity of Lyapunov exponents with respect to the cocycle and the invariant probability for random products of $SL(2,\R)$ matrices in the Bernoulli setting. In higher dimensions,
continuous dependence with respect to $A$ of all Lyapunov exponents for i.i.d. random products of matrices in $GL(d,\R)$ was proved by Avila et al. \cite{ave}.  If the base dynamics is a subshift
of finite type or, more generally, a hyperbolic set, then Backes-Brown-Butler \cite{bbb} proved
that the Lyapunov exponents are always continuous among H\"older continuous fiber-bunched
$SL(2,\R)$-cocycles.

If $A\in C^r(X,SL(2,\R))$, it is known that $L(A)$ is upper semicontinuous; thus, it is continuous at generic $A$. Especially, it is continuous at A with $L(A)=0$ and at uniformly hyperbolic cocycles. The most interesting issue is the continuity of $L(A)$ at the nonuniformly hyperbolic cocycles, which is found to depend on the class of cocycles under consideration including its topology. LE was proved to be discontinuous at any nonuniformly hyperbolic cocycles in $C^0$-topology
by Furman \cite{furman} (Continuity at uniformly hyperbolic cocycles is well-known). Motivated by
Mane \cite{mane1,mane2}, Bochi \cite{Bochi} further proved that any nonuniformly hyperbolic $SL(2,\R)$-cocycle over a fixed ergodic system on a compact space, can be arbitrarily approximated
by cocycles with zero LE in the $C^0$-topology.

In this paper, we are interested in the quasiperiodic cocycles. The base system is a rotation on the torus  in this case, things are very complicated: it will depend on the smoothness of $A$ in a very sensitive way. If the cocycle is analytic,  the  H\"{o}lder continuity of the Lyapunov exponent in the positive Lyapunov exponent regime was proved by Goldstein and Schlag \cite{gs1} assuming that $\alpha$ is  strong Diophantine. Similar results were proved
in \cite{bgs} by Bourgain, Goldstein, and Schlag when the underlying dynamics is a shift
or skew-shift of a higher-dimensional torus. For more results of this favor, here is a partial list \cite{b1,gs1,gs2,FT,hz,gyzh,LWY,PV,wz1,YZ,xgw,s,dk1,dk2}.  Later, it was proved by Bourgain-Jitomirskaya  in \cite{bj} that the LE is joint continuous for $SL(2,\R)$ cocycles, in frequency and cocycle map, at any irrational frequencies.  Jitomirskaya-Koslover-Schulteis \cite{jks} got the continuity of LE with respect to potentials for a class of
analytic quasiperiodic $M(2,\C)$-cocycles.  Bourgain \cite{bourgain2} extended the results in \cite{bj} to multi-frequency case.
Jitomirskaya-Marx \cite{jm} extended the results in \cite{bj} to all (including singular) $M(2,\C)$-cocycles. More recently, continuity of the Lyapunov exponents for one-frequency analytic $M(m,\C)$ cocycles was given by Avila-Jitomirskaya-Sadel \cite{ajs}. Weak H\"older continuity of the Lyapunov exponents for multi-frequency $GL(m, \C)$-cocycles, $m \geq 2,$ was recently obtained by Schlag \cite{s} and
Duarte-Klein \cite{dk1}. For the lower regularity case, Klein \cite{klein} proved that for Schr\"{o}dinger operators with potentials in a Gevrey class $G^s$ with $1\leq s<2$, the LE  is  weak H\"{o}lder continuous on any compact interval of the energy provided that  the frequency is strong Diophantine and the LE is large than $0$. While if we further lower the regularity of the potential, Wang-You \cite{wangyou} constructed examples to show
that the LE of quasiperiodic Schr\"odinger cocycles can be discontinuous with respect to the potential even in the $C^{\infty}$-topology. Jitomirskaya-Marx \cite{jm} obtained similar results in the complex category $M(2,\C)$ by the tools of harmonic analysis. Recently,
Wang-You \cite{wy2}  improved the result in \cite{wangyou} by showing that in $C^r$-topology, $1\leq r\leq +\infty,$ there exists Schr\"odinger cocycles with a positive LE that can be approximated by
ones with zero LE. For other results about results on discontinuity of LE, one can see \cite{DGK,viana2}.

\subsection{Outline of the proof and the structure of this paper}
The main results of this paper are based on several improvements of the results in \cite{wangyou} where the authors constructed examples of discontinuity of LE in $C^\infty$-topology. We first give a quick review of the main ideas. We construct $D_s\in G^s(\mathbb{S}^1,SL(2,\R))$ $(s>2)$ as the limit of a sequence of cocycles $\{A_{n},n=N,N+1,...\}$ in $G^{s}(\mathbb{S}^1,SL(2,\R))$ $(s>2)$, the sequence  $\{A_{n},n=N,N+1,...\}$ possesses some kind of finite hyperbolic property, that is, $\|A_n^{r_n^+}(x)\|\approx \lambda^{r_n^+}$ for most $x\in\mathbb{S}^1$ and $\lambda\gg1$ with $r_n^+\rightarrow\infty$ as $n\rightarrow\infty$, which gives a lower bound estimate $(1-\e)\ln\lambda$
of the Lyapunov exponent of the limit cocycle $(\alpha,D_s)$. Then we modify $\{A_{n},n=N,N+1,...\}$, and construct another sequence of cocycles $\{\widetilde{A}_{n},n=N,N+1,...\}$ with some kind of degenerate property such that $\widetilde{A}_n\rightarrow D_s$ in $G^s$-topology as $n\rightarrow\infty$. Moreover, for each $n$, the Lyapunov
exponent of $(\alpha,\widetilde{A}_n)$ is less than $(1-\delta)\ln\lambda$ with $\delta\gg\e$,
which implies the discontinuity of the Lyapunov exponent at $D_s$.

Compared to \cite{wangyou}, the main technical improvements of the present paper are the following two aspects:
\begin{enumerate}
\item Since we need to construct examples in Gevrey space, we need explicit examples of Gevrey functions. We find the $C^\infty$-bump functions are all Gevrey functions based on an optimal estimate of the upper bound of its derivatives. Surprisingly, this easy but important observation makes it possible for us to construct a counterexample in Gevrey space.
\item Another technical difficulty (the most difficult part) is to prove the sequences $\{A_n\}_{n=N}^\infty$ and $\{\widetilde{A}_n\}_{n=N}^\infty$ converge in $G^s$-topology ($s>2$). It is much more difficult than to prove the convergence in $C^\infty$-topology since one needs very delicate control of the derivatives, and it is out of reach by the methods in \cite{wangyou}.  We overcome this difficulty by developing a $G^s$ version of the concatenation of finitely many hyperbolic matrices, i.e., Lemma \ref{Gs1} and Lemma \ref{2} in our paper. Our new Lemmas enable us to not only greatly simplify the proofs in \cite{wangyou}, but also optimize almost all estimates in \cite{wangyou}. Our construction is optimal since \cite{cgyz} has shown that for the case of $s<2$, the Lyapunov exponent is continuous.
\end{enumerate}
A key technique in the construction of $A_n(x)$ comes from Young \cite{young}, which was derived from Benedicks and Carleson \cite{bc}.  Based on this technique, Wang-Zhang \cite{wz1} developed a new iteration scheme to prove LDT for Schr\"odinger cocycles with a class of finitely differential potential. They proved that for $C^2$ cos-like (Morse) potential with a large coupling, LE is weak-H\"older continuous. In this aspect, we also give an improvement of the non-resonance lemma proved in \cite{wz1}, which plays an important role in our proof. For more applications of Benedicks-Carleson-Young's method to quasiperiodic Schr\"odinger operators, we refer readers to \cite{bk1,bk2,wz2}.

The structure of this paper is as follows. In Section 2, we give some basic concepts and preparations for $SL(2,\R)$ matrices and Gevrey functions. The main idea of the proof will be sketched in Section 3. In Section 4, we give the details of the construction, which is the key part of this paper.  Finally, We give the proof of some basic properties of Gevrey functions in Section 5.
\section{Preparations and some technical lemmas}
For $\theta \in \mathbb{S}^1,$ let
$$
R_{\theta}=\begin{pmatrix}
\cos{\theta}&-\sin{\theta}\\
\sin{\theta}&\cos{\theta}
\end{pmatrix}
\in SO(2,\R).
$$
Define the map
$$s:SL(2,\R)\rightarrow\R{\mathbb{P}}^1=\R/(\pi\Z)$$
so that $s(A)$ is the most contraction direction of $A\in SL(2,\R).$ That is, for a unit vector $\hat{s}(A)\in s(A)$, it holds that $\|A\cdot\hat{s}(A)\|=\|A\|^{-1}.$ Abusing the notation a little, let
$$u:SL(2,\R)\rightarrow \R{\mathbb{P}}^1=\R/(\pi\Z)$$
be determined by $u(A)=s(A^{-1})$ and $\hat{u}(A)\in u(A)$. Then for $A\in SL(2,R)$, it is clear that
$$
A=R_u\cdot\begin{pmatrix}
\|A\|&0\\
0&\|A\|^{-1}
\end{pmatrix}
\cdot R_{\frac{\pi}{2}-s},
$$
where $s,u\in[0,\pi)$ are angles corresponding to the directions $s(A),u(A)\in \R/(\pi\Z).$

\subsection{Hyperbolic sequences of $SL(2,\R)$-matrices}

For a sequence of matrices $\{\cdots A_{-1},A_0,A_1,\cdots\}$, we denote
$$
A^n=A_{n-1}\cdots A_1A_0
$$
and
$$
A^{-n}=A^{-1}_{-n}\cdots A_{-1}^{-1}.
$$
\begin{Definition}
For any $1<\mu\leq \lambda$, we say the block of matrices $\{A_0,A_1,\cdots,A_{n-1}\}$ is $\mu$-hyperbolic if
\begin{enumerate}[(1)]
\item $\|A_i\|\leq \lambda,\ \ \forall i$,
\item $\|A^i\|\geq \mu^{i(1-\e)},\ \ \forall i$,
\end{enumerate}
and (1) and (2) hold if $A_0,\cdots,A_{n-1}$ are replaced by  $A^{-1}_{n-1},\cdots,A^{-1}_{0}$.
\end{Definition}
The following lemma is due to Young \cite{young} which tells us when the
concatenation of two hyperbolic blocks is still a hyperbolic block.
\begin{Lemma}[Lemma 5 of \cite{young}]\label{2.1}
Suppose that $C$ satisfies $\|C\|\geq \mu^m$ with $\mu\gg 1$. Assume that $\{A_0,A_1,\cdots,A_{n-1}\}$ is a $\mu$-hyperbolic sequence, and assume that $\angle(s(C^{-1}),s(A^n))=2\theta\ll 1$. Then $\|A^n\cdot C\|\geq \mu^{(m+n)(1-\e)}\cdot\theta$.
\end{Lemma}
\subsection{The Gevrey functions}
\subsubsection{Basic properties}
In the following, $s\geq 1, K>0$ will always be some fixed constants. We give some basic properties on the product and composition of Gevrey functions whose proofs will be given in Section 5. Abusing the notations a little bit, for any $f\in G^{s,K}(I)$, we denote
$$
|f|_{s,K}=|f|_{G^{s,K}(I)}:=\frac{4\pi^2}{3}\sup\limits_{k}\frac{(1+|k|)^2}{K^k(k!)^s}|\partial^k f|_{C^0(I)}.
$$
\begin{Proposition}\label{progev}
Assume $f,g\in G^{s,K}(I)$ and  $\e>0$ is sufficiently small, we have
\begin{enumerate}[$(1)$]
\item{ $|fg|_{s,K}\leq |f|_{s,K}|g|_{s,K}$.}
\item{ For any $\e>0$, $\left|\partial f\right|_{s,(1+\e^{\frac{1}{s}})K}\leq \frac{K}{\e}.$}
\item{ Assume $|f-1|_{s,K}\leq \e$, then
$$
\left|\frac{1}{f}-1\right|_{s,(1+\e^{\frac{1}{s+8}})K}\leq \e^{\frac{1}{12}}.
$$}
\item{ Assume $|f-1|_{s,K}\leq \e$, then
$$
\left|\sqrt{f}-1\right|_{s,(1+\e^{\frac{1}{s+16}})K}\leq \e^{\frac{1}{12}}.
$$}
\item{ Assume $|f|_{s,K}\leq \e$, then $\arcsin(f)\in G^{s,4K}(I)$ and
$$
\left|\sin f\right|_{s,(1+\e^{\frac{1}{s+8}})K},\ \ \left|\cos f-1\right|_{s,(1+\e^{\frac{1}{s+8}})K}\leq \e^{\frac{1}{12}}.
$$}
\end{enumerate}
\end{Proposition}

\subsubsection{Explicit examples}
Given a $C^\infty$-bump function, an interesting question is to investigate the decay rate of its Fourier coefficients (equivalently, the growth rate of its derivatives). In this part, we investigate the Gevrey exponent of various $C^\infty$-bump functions, the proofs will be also postponed to Section 5. We remark that our estimates of the upper bound of the derivatives of $C^\infty$ bump functions are even optimal, we refer readers to \cite{steven} for more details.
\begin{Lemma}\label{Gevrey-function}
Assume that $0<\nu<\infty$ and
\begin{equation*}
f(x)=
\begin{cases}
e^{-\frac{1}{|x|^\nu}}& x\neq 0\\
0& x=0,
\end{cases}
\end{equation*}
then there is some $C>0$ such that
$$
|f^{(n)}(x)|\leq \begin{cases}C^n e^{-\frac{1}{2|x|^\nu}}(n!)^{1+\frac{1}{\nu}} &x\neq 0\\
0& x=0
\end{cases}
,\ \  \forall n\in\N.
$$
As a corollary, $f\in G^{1+\frac{1}{\nu}}(\R)$.
\end{Lemma}

\begin{Corollary}\label{inverse-Gevrey-function}
Assume that $0<\nu<\infty$ and
\begin{equation*}
f(x)=e^{\frac{1}{|x|^\nu}},\ \  x\neq 0.
\end{equation*}
Then there is some $C>0$ such that for any $x\neq 0$, we have
$$
|f^{(n)}(x)|\leq C^ne^{\frac{2}{|x|^\nu}}(n!)^{1+\frac{1}{\nu}}, \ \ \forall n\in\N.
$$
\end{Corollary}
\begin{Corollary}\label{phi_0}
We define a $2\pi$-periodic function as follows
$$
g(x)=\begin{cases}ce^{-\left(\frac{1}{(x-c_1-k\pi)^\nu}+\frac{1}{(c_1+(k+1)\pi-x)^\nu}\right)}& x\in (c_1+k\pi,c_1+(k+1)\pi)\\
0& x\in \{c_1+k\pi, c_1+(k+1)\pi\}
\end{cases},
$$
then there is some $C>0$ such that for any $x\in \mathbb{S}^1$, we have
$$
|g^{(n)}(x)|\leq C^n e^{-\left(\frac{1}{2|x-c_1|^\nu}+\frac{1}{2|x-c_1-\pi|^\nu}\right)}(n!)^{1+\frac{1}{\nu}},\ \  \forall n\in\N.
$$
As a corollary, $g\in G^{1+\frac{1}{\nu}}(\mathbb{S}^1)$.
\end{Corollary}

Let $I_{n,1}=\left[c_1-\frac{1}{q_n^\beta},c_1+\frac{1}{q_n^\beta}\right]$, $I_{n,2}=\left[c_2-\frac{1}{q_n^\beta},c_2+\frac{1}{q_n^\beta}\right]$ and $I_n=I_{n,1}\bigcup I_{n,2}$.
\begin{Lemma}\label{fn}
Assume $0<\nu<1$, $\beta>1$ and $\delta>0$ satisfy $0<\frac{\beta}{\frac{1}{\nu}-\delta}<1$. For any $n\geq N$, there exist an absolute constant $C$ and a $2\pi$-periodic function $f_n\in G^{1+\frac{1}{\nu},C}(\mathbb{S}^1)$ such that
$$
f_n(x)\left\{
\begin{aligned}
&=1&x\in \frac{I_n}{10}\\
&\in (0,1] &x\in I_n\backslash\frac{I_n}{10}\\
&=0& x\in \mathbb{S}^1\backslash I_n
\end{aligned}
\right.
$$
and
$$
|f_n|_{1+\frac{1}{\nu},C}\leq (Cq^\beta_n)^{q_n^{\frac{\nu\beta}{1-\delta\nu}}}.
$$
\end{Lemma}

\section{Proof of Theorem \ref{main1}}
We first introduce some notations. Let $p_n/q_n$ be the continued fraction expansion of $\alpha$. The general settings of the cocycle $(\alpha,A)$ are
\begin{itemize}
\item $q_{n+1}\leq Mq_n$, $n\in\N$,
\item $A\in G^{1+\frac{1}{\nu}}(\mathbb{S}^1)$ with $0<\nu<1$.
\end{itemize}

Let $M,N>0$ be sufficiently large such that
$$
\e=M^{-100}\ll\delta=\frac{1}{4}M^{-20},\ \ \lambda=e^{q_N ^{q_N}}\gg 1,\ \ 1<\beta<\frac{1}{\nu}.
$$
Denote by $\gamma=\nu\beta$. For $n\geq N$, we inductively define
$$
\ln\lambda_{n+1}=\ln\lambda_n-10^4 q_{n+1}^{\gamma-1}, \ \ \lambda_N=\lambda^{1-\e}.
$$
$$
\ln\widetilde{\lambda_{n+1}}=\ln\widetilde{\lambda_n}+10^4 q_{n+1}^{\gamma-1},\ \ \widetilde{\lambda_N}=\lambda^{1+\e}.
$$
Choose $q_N$ sufficiently large such that $\sum\limits_{i=N+1}^{\infty}10^4 q_i^{\gamma-1}<\e$, then $\lambda_\infty \geq \lambda^{1-2\e}$ and $\widetilde{\lambda_\infty} \leq \lambda^{1+2\e}$.

We define
\begin{itemize}
\item The critical set: $C_0=\{c_1,c_2\}$ where $c_1\in [0,\pi)$ and $c_2=c_1+\pi$.
\item The critical interval: $I_{n,1}=\left[c_1-\frac{1}{q_n^\beta},c_1+\frac{1}{q_n^\beta}\right]$, $I_{n,2}=\left[c_2-\frac{1}{q_n^\beta},c_2+\frac{1}{q_n^\beta}\right]$ and $I_n=I_{n,1}\bigcup I_{n,2}$.
\item The first return time: For $x\in I_n$, we denote the smallest positive integer $i$ with $T^ix\in I_n$ (respectively $T^{-i}x\in I_n$) by $r_n^+(x)$ (respectively $r_n^-(x)$), and define $r_n^{\pm}=\min_{x\in I_n}r_n^{\pm}(x)$. Obviously, $r_n^{\pm}\geq\frac{q_n}{2}$.
\item The sample function: The $2\pi$-periodic smooth function $\phi_0$ is defined as
$$
\sin \phi_0(x)=ce^{-\left(\frac{1}{(x-c_1-k\pi)^\nu}+\frac{1}{(c_1+(k+1)\pi-x)^\nu}\right)},\ \ x\in [c_1+k\pi,c_1+(k+1)\pi),\ \ k\in\Z,
$$ where $c$ is sufficiently small.
\end{itemize}
\begin{Remark}
To ensure $r_n^{\pm}\geq \frac{q_n}{2}$, we must require $\beta>1$. To ensure $\sum\limits_{i=N+1}^{\infty}10^4 q_i^{\gamma-1}<\e$, we must require $\beta\nu<1$. Thus our construction is possible only if $\nu<1$. Indeed, it is essential since if $\nu>1$, the LE is continuous \cite{cgyz,klein}.
\end{Remark}
\begin{Remark}
By (5) in Proposition \ref{progev} and Corollary \ref{phi_0}, we have $\phi_0\in G^{1+\frac{1}{\nu}}(\mathbb{S}^1)$.
\end{Remark}
For $C\geq 1$, we denote by $\frac{I_{n,i}}{C}$ the sets $\left[c_i-\frac{1}{Cq_n^\beta},c_i+\frac{1}{Cq_n^\beta}\right]$, $i=1,2$ and by $\frac{I_n}{C}$ the set $\frac{I_{n,1}}{C}\bigcup \frac{I_{n,2}}{C}$. Let $\Lambda=\begin{pmatrix}\lambda&0\\ 0&\lambda^{-1}\end{pmatrix}.$ Theorem \ref{main1} follows from the following two propositions whose proofs will be given in Section 4.
\begin{Proposition}\label{p1}
There exist  functions  ${\phi}_n(x)$ on $\mathbb{S}^1$ $(n=N,N+1,\cdots)$ such that\\

\begin{enumerate}[$(1)_n$]
\item $|\phi_n-\phi_{n-1}|_{1+\frac{1}{\nu},K}\leq \lambda_n^{-\frac{q_{n-1}}{100}}$ for some $K>0$, if $n>N$,\\

\item $A_n(x),A_n(Tx),\cdots,A_n(T^{r_n^+(x)-1}(x))$ is $\lambda_n$-hyperbolic  for $x\in I_n$ where  ${A}_n(x)=\Lambda R_{\frac{\pi}{2}-{\phi}_n(x)}$,\\

\item  It holds

\ \ \ \ $(a)_n$\ \ $s_n(x)-s_n^\prime(x)=\phi_0(x)\ \  x\in \frac{I_n}{10},$

\ \  \ \ $(b)_n$\ \ $|s_n(x)-s_n^\prime(x)|\geq \frac{1}{2}|\phi_0(x)|\geq ce^{-10^\nu q_n^{\gamma}}, \ \ x\in I_n\backslash \frac{I_n}{10},$\\
\ \  \ \ \ \ \ \ where $s_n(x)=\overline{s(A_n^{r_n^+}(x))}, s'_n(x)=\overline{s(A_n^{-r_n^-}(x))}$,\\
\item  It holds
$$
\|A_{r_n^{\pm}}\|_{G^{1+\frac{1}{\nu},K}(I_n)}\leq \widetilde{\lambda_n}^{r_n^{\pm}},\ \ \left|\frac{1}{\|A_{r_n^{\pm}}\|}\right|_{G^{1+\frac{1}{\nu},K}(I_n)}\leq \lambda_n^{-r_n^{\pm}}.
$$
\end{enumerate}
\end{Proposition}

\begin{Proposition}\label{p2}
There exist  functions  $\tilde{\phi}_n(x)$ on $\mathbb{S}^1$ $(n=N,N+1,\cdots)$ such that \\
\begin{enumerate}[$(1)_n$]
\item  $|\phi_n(x)-\tilde{\phi}_{n}(x)|_{1+\frac{1}{\nu},K}\leq Cq_n^{-2}$ for some $K>0$, if $n>N$,\\

\item $\widetilde{A}_n(x),\widetilde{A}_n(Tx),\cdots,\widetilde{A}_n(T^{r_n^+(x)-1}(x))$ is $\lambda_n$-hyperbolic for each $x\in I_n$ where $\widetilde{A}_n(x)=\Lambda R_{\frac{\pi}{2}-\tilde{\phi}_n(x)}$,\\

\item It holds
\begin{align*}
\tilde{s}_n(x)=\tilde{s}_n^\prime(x),\ \  x\in \frac{I_n}{10},
\end{align*}
 where $\tilde{s}_n(x)=\overline{s(\widetilde{A}_n^{r_n^+}(x))}$, $\tilde{s}_n^\prime(x)=\overline{s(\widetilde{A}_n^{-r_n^-}(x))}$.
\end{enumerate}
\end{Proposition}
\noindent
{\bf Proof of Theorem \ref{main1}:} By (1)'s in Proposition \ref{p1} and Proposition \ref{p2}, there exists $D_{\nu}\in G^{1+\frac{1}{\nu},K}$ such that $\|A_n-D_\nu\|_{1+\frac{1}{\nu},K}\rightarrow 0$ and $\|\widetilde{A}_n-D_\nu\|_{1+\frac{1}{\nu},K}\rightarrow 0$. Then Theorem \ref{main1} is a direct conclusion of the followings:\\
(a) $L(D_\nu)\geq (1-4\e)\ln\lambda$,\\
(b) $L(\widetilde{A}_n)\leq (1-\delta)\ln\lambda,\ \ \forall n> N$.\\
\\
Step 1: Proof of (a). We say $x\in\mathbb{S}^1$ is nonresonant for $A_n(x)$ if
\begin{equation}\label{nonresonant}
\begin{cases}
dist (T^ix,C_0)>\frac{1}{q_N^\beta},& 0\leq i<q_N,\\
dist (T^ix,C_0)>\frac{1}{q_k^\beta},& q_{k-1}\leq i<q_k, N<k\leq n.
\end{cases}
\end{equation}
The Lebesgue measure of the set of all nonresonant points  $x\in\mathbb{S}^1$ is at least $2\pi(1-\sum_{N\leq k<n}1/q^{\beta-1}_k)$, which is larger than $2\pi(1-\e/2C\pi)$ for $N\gg1$. For any $x$ satisfying the nonresonant property \eqref{nonresonant}, let $j_0$ be the first time such that $T^jx\in I_N$ and let $n_0$ be such that $T^{j_0}x\in I_{n_0}\backslash I_{n_0+1}$. In general, let $j_i$ and $n_i$ be defined so that $T^{j_i}x\in I_{n_i}\backslash I_{n_i+1}$ and let $T^{j_{i+1}}x$ be the next return of $T^{j_i}x$ to $I_{n_i}$. It is obvious that $j_{i+1}-j_i\geq \frac{q_{n_i}}{2}$. By condition (2), we have $\{A_{n}(T^{j_i}x),\cdots,A_{n}(T^{j_{i+1}-1}x)\}$ is $\lambda_\infty$-hyperbolic (See also \cite{young} for similar arguments).

Since $T^{j_i}x\notin I_{n_i+1}$, by $(3)_n$ of Proposition \ref{p1}  and the definition of $\phi_0$, we have
$$
\angle(s_{n}(T^{j_i}x),s'_{n}(T^{j_i}x))\geq ce^{-10^\nu q_{n_i+1}^{\gamma}}.
$$
On the other hand, it holds that
$$
\angle\left(s(A_{n}^{-j_i}(T^{j_i}x)),s(A_{n}^{j_{i+1}-j_i}(T^{j_i}x))\right)>\frac{1}{2}\angle\left(s_{n}(T^{j_i}x),s'_{n}(T^{j_i}x)\right)\geq ce^{-10^\nu q_{n_i+1}^{\gamma}}.
$$
By Lemma \ref{2.1}, we have
\begin{align*}
\|A_n^{j_{i+1}}(x)\|&\geq \|A_n^{j_i}(x)\|\cdot\|A_n^{j_{i+1}-j_i}(T^{j_i}x)\|\cdot \angle(s(A_{n}^{-j_i}(T^{j_i}x)),s(A_{n}^{j_{i+1}-j_i}(T^{j_i}x)))\\
&\geq \|A_{n}^{j_i}(x)\|\cdot\lambda_\infty^{(j_{i+1}-j_i)(1-\e)}ce^{-10^\nu q_{n_i+1}^{\gamma}}.
\end{align*}
Inductively
$$
\|A_n^{j_s}(x)\|\geq \|A_n^{j_0}(x)\|\cdot\lambda_\infty^{j_s-j_0}\cdot\prod_{i=0}^{s-1}ce^{-10^\nu q_{n_i+1}^{\gamma}}.
$$
Notice that $j_s-j_0=\sum\limits_{i=1}^{s}j_{i}-j_{i-1}\geq \sum\limits_{i=0}^{s-1}\frac{q_{n_i}}{2}$, thus if $s$ is sufficiently large, we have
$$
j_0\leq q_N^C\leq \frac{\e}{4}j_s,\ \ \prod_{i=0}^{s-1}ce^{-10^\nu q^\gamma_{n_i+1}}\geq \lambda^{-\frac{\e}{4} (j_s-j_0)}.
$$
It follows that
$$
\|A_n^{j_s}(x)\|\geq \lambda_\infty^{(1-\frac{3\e}{2})j_s} \geq \lambda^{(1-2\e)j_s} _{\infty}.
$$

Now we are ready to prove the main result. From the subadditivity of the cocycle, the finite Lyapunov exponent of a cocycle converges to the Lyapunov exponent. Thus there exists a large $s_0\geq N_0$ such that
$$
\left|\frac{1}{j_{s_0}}\int_{\T}\ln\|D^{j_{s_0}}_\nu(x)\|dx-L(\alpha,D_\nu)\right|\leq \e.
$$
By $(1)_n$ of Proposition \ref{p1}, there exists $N_1>N_0$, such that for any $n>N_1$, it holds that
$$
\left|\frac{1}{j_{s_0}}\int_{\T}\ln\|D^{j_{s_0}}_\nu(x)\|dx-\frac{1}{j_{s_0}}\int_{\T}\ln\|A_n^{j_{s_0}}(x)\|dx\right|\leq \e.
$$
On the other hand
$$
\frac{1}{j_{s_0}}\int_{\T}\ln\|A_n^{j_{s_0}}(x)\|dx\geq (1-2\e)\ln\lambda_\infty-C\frac{\e}{C}\geq (1-3\e)\ln\lambda.
$$
Thus we finish the proof of (a).\\
\\
\noindent
Step 2: Proof of (b). The following two lemmas have been proved in \cite{wangyou}.
\begin{Lemma}[Lemma 4.1 of \cite{wangyou}]\label{u1}
Suppose $A$ and $B$ are two hyperbolic matrices such that $\|A\|=\lambda_1^m$ and $\|B\|=\lambda_2^n$ with $m,n>0$ and $\lambda_1,\lambda_2\gg 1$. If $A(s(A)) || u(B)$, then $\|BA\|\leq 2\max\{\lambda_1^m\lambda_2^{-n},\lambda_2^n\lambda_1^{-m}\}$.
\end{Lemma}
\begin{Lemma}[Corollary 4.1 of \cite{wangyou}]\label{u2}
Let $\min r_n=\min_{x\in I_n}\min\{i>0| T^ix(\mod 2\pi)\in I_n\}$, and let $\max r_n=\max_{x\in \frac{1}{10}I_n}\min\{i>0| T^ix(\mod 2\pi)\in \frac{I_n}{10}\}$, Then $M^{-20}\leq \frac{\min r_n}{\max r_n}\leq 1$.
\end{Lemma}
Let $\cdots<n_{j-1}<n_j<n_{j+1}<\cdots$ be the returning times of $x\in \frac{I_n}{10}$ to $\frac{I_n}{10}$. Moreover, we let $n_{j+}$ be the first returning time of $x\in I_n$ to $I_n$ after $n_j$. Similarly, we denote by $n_{j-}$  the last returning time of $x\in I_n$ to $I_n$ before $n_j$. Obviously, it holds that $n_{j-1}\leq n_{j-}<n_j$ and $n_j<n_{j+}\leq n_{j+1}$. By Lemma \ref{u2}, we have
$$
n_{j_+}-n_j, n_j-n_{j_-}\leq (1-\frac{1}{2}M^{-20})(n_{j_+}-n_{j_-}).
$$

Since $T^{n_j}x\in I_n$, by  Proposition \ref{p2} and Lemma \ref{u1}, we have
\begin{align*}
\left\|\widetilde{A}_n(T^{n_{j_+}}x)\cdots\widetilde{A}_n(T^{n_{j}}x)\cdots\widetilde{A}_n(T^{n_{j_-}}x)\right\|&\leq2 \max\left\{\widetilde{A}_n^{n_{j_+}-n_j}(T^{n_jx}),\widetilde{A}_n^{n_{j}-n_{j_-}}(T^{n_{j_-}x})\right\}\\
&\leq 2\lambda^{\max\{n_{j_+}-n_j,n_j-n_{j_-}\}}\leq \lambda^{(1-\frac{1}{2}M^{-20})(n_{j_+}-n_{j_-})}.
\end{align*}
It follows that
$$
\left\|\widetilde{A}_n(T^{n_{j+1}}x)\cdots\widetilde{A}_n(T^{n_{j-1}}x)\right\|\leq \lambda^{n_{j+1}-n_{j-1}-\frac{1}{2}M^{-20}(n_{j_+}-n_{j_-})}\leq  \lambda^{(n_{j+1}-n_{j-1})(1-\frac{1}{4}M^{-40})}.
$$
Thus for any even $k$,
$$
\left\|\widetilde{A}_n(T^{n_{k}}x)\cdots\widetilde{A}_n(x)\right\|\leq \lambda^{(1-\frac{1}{4}M^{-40})\sum_{j=0}^{\frac{k}{2}-1}(n_{2j+2}-n_{2j})} \leq \lambda^{n_k(1-\frac{1}{4}M^{-40})},
$$
which implies that $L(\widetilde{A}_n)\leq (1-\delta)\ln\lambda$.

\section{Proof of Proposition \ref{p1} and \ref{p2}}
In this section, we aim to prove Proposition \ref{p1} and Proposition \ref{p2} which are the main technical parts of this paper. The proof is split into the following three subsections.
\subsection{Key Lemmas}
Let $2<s<\infty$, $0<\gamma<1$, $K>0$, $\alpha\in\R\backslash\Q$ be bounded,  $p_n/q_n$ be the continued fraction expansion of $\alpha$ and $\lambda>1$ be sufficiently large. Without loss of generality, we may assume $2<s<3$ since we are mainly concerned with the case $s$ is sufficiently close to $2$. Recall that for any $f\in C^\infty(I)$, we denote
$$
|f|_{s,K}:=\frac{4\pi^2}{3}\sup\limits_{k}\frac{(1+|k|)^2}{K^k(k!)^s}|\partial^k f|_{C^0}(I),
$$
the Gevrey norm of $f$ restricting to $I$.

In the following, we prove a Gevrey version of the concatenation of hyperbolic matrices. It greatly simplifies and improves the proofs in \cite{wangyou}.
\begin{Lemma}\label{Gs1}
Let
$$
E(x)=\begin{pmatrix}
e_2(x)&0\\
0&e_2^{-1}(x)
\end{pmatrix}R_{\theta(x)}\begin{pmatrix}
e_1(x)&0\\
0&e_1^{-1}(x)
\end{pmatrix},
$$ where  $e_1, e_2,\theta \in G^{s,K}(I)$ satisfying
\begin{equation}\label{lowcond}
\inf\limits_{x\in I}\left|\theta(x)-\frac{\pi}{2}\right|\geq ce^{-q_n^\gamma}\gg \min\left\{\inf\limits_{x\in I} e_1(x),\inf\limits_{x\in I} e_2(x)\right\}^{-\frac{1}{100}},
\end{equation}
\begin{equation}\label{lowcond1}
\left|\frac{1}{\cos \theta}\right|_{s,K}, \left|\tan \theta\right|_{s,K}\leq Ce^{q_n^\gamma}, \ \ \left|\cos\theta\right|_{s,K},\left|\cot\theta\right|_{s,K}\leq C,
\end{equation}
\begin{equation}\label{lowcond2}
|e^{-1}_i|_{s,K}\leq C\lambda^{-\frac{1}{3}q_{n-1}},\ \ i=1,2.
\end{equation}
Then, for $e(x):=\|E(x)\|$, it holds that
\begin{equation}\label{r3}
\inf\limits_{x\in I} e(x)\geq c\inf\limits_{x\in I} e_1(x)\cdot \inf\limits_{x\in I} e_2(x)\cdot e^{-q_n^\gamma},
\end{equation}
\begin{equation}\label{r4}
|e|_{s,(1+\eta)K}\leq C^2|e_1|_{s,K}|e_2|_{s,K},
\end{equation}
\begin{equation}\label{r1}
|e^{-1}|_{s,(1+\eta)K}\leq C^2|e^{-1}_1|_{s,K}|e^{-1}_2|_{s,K}e^{q_n^\gamma}.
\end{equation}
Let $s(x)=s(E(x))$ and $u(x)=u(E(x))$, we further have
\begin{equation}\label{r2}
\left|\frac{\pi}{2}-s\right|_{s,(1+\eta)K}\leq |e_1^{-1}|^{3}_{s,K}|e_1|^2_{s,K},\ \ \left|u\right|_{s,(1+\eta)K}\leq |e_2^{-1}|^{3}_{s,K}|e_2|^2_{s,K},
\end{equation}
where $\eta=\lambda^{-\frac{1}{4000}q_{n-1}}$.
\end{Lemma}
\begin{pf}
For simplicity, let us omit the dependence on $x$ in the following computation. Direct computations show that
$$
E^tE=\begin{pmatrix}
e_1^2e_2^2\cos^2\theta+e_1^2e_2^{-2}\sin^2\theta& (e_2^{-2}-e_2^2)\sin\theta\cos\theta\\
(e_2^{-2}-e_2^2)\sin\theta\cos\theta&e_1^{-2}e_2^{-2}\cos^2\theta+e_2^2e_1^{-2}\sin^2\theta
\end{pmatrix}.
$$
It is obvious that
\begin{align*}
e^2+e^{-2}=&e_1^2e_2^2\cos^2\theta+e_1^2e_2^{-2}\sin^2\theta+e_1^{-2}e_2^{-2}\cos^2\theta+e_2^2e_1^{-2}\sin^2\theta\\
=&e_1^2e_2^2\cos^2\theta\left(1+e_2^{-4}\tan^2\theta+e_1^{-4}e_2^{-4}+e_1^{-4}\tan^2\theta\right):=b.
\end{align*}
Thus
\begin{equation}\label{re_e}
e=\sqrt{\frac{b+\sqrt{b^2-4}}{2}}=\sqrt{b}\sqrt{\frac{1+\sqrt{1-4b^{-2}}}{2}},
\end{equation}
\begin{equation}\label{re_e1}
e^{-1}=\sqrt{\frac{2}{b+\sqrt{b^2-4}}}=\sqrt{\frac{1}{b}}\sqrt{\frac{2}{1+\sqrt{1-4b^{-2}}}}.
\end{equation}
By \eqref{lowcond} and \eqref{re_e}, we have
\begin{equation}\label{ne6}
\inf\limits_{x\in I} e(x)\geq c\inf\limits_{x\in I} e_1(x)\cdot \inf\limits_{x\in I} e_2(x)\cdot e^{-q_n^\gamma}.
\end{equation}
By (1) in Proposition \ref{progev}, \eqref{lowcond2} and \eqref{lowcond1}, we have
\begin{align}\label{ne1}
\nonumber &\left|e_2^{-4}\tan^2\theta+e_1^{-4}e_2^{-4}+e_1^{-4}\tan^2\theta\right|_{s,K}\\ \nonumber
\leq& |e_2^{-1}|_{s,K}^4|\tan\theta|_{s,K}^2+|e_1^{-1}|_{s,K}^4|e_2^{-1}|_{s,K}^4+|e_1^{-1}|_{s,K}^4|\tan\theta|_{s,K}^2\\
\leq &C\lambda^{-\frac{4}{3}q_{n-1}}e^{2q_n^\gamma}+C\lambda^{-\frac{8}{3}q_{n-1}}+C\lambda^{-\frac{4}{3}q_{n-1}}e^{2q_n^\gamma}\leq \lambda^{-q_{n-1}}.
\end{align}
The last inequality holds since $\alpha$ is bounded and $\lambda$ is sufficiently large.

Let $\eta_1=\lambda^{-\frac{1}{20}q_{n-1}}$ and $\eta_2=\lambda^{-\frac{1}{240}q_{n-1}}$, by \eqref{ne1} and (3)-(4) in Proposition \ref{progev}, we have
\begin{align}\label{ne2}
\left|\frac{1}{1+e_2^{-4}\tan^2\theta+e_1^{-4}e_2^{-4}+e_1^{-4}\tan^2\theta}-1\right|_{s,(1+\eta_1)K}\leq \lambda^{-\frac{1}{12}q_{n-1}}.
\end{align}
\begin{align}\label{ne9}
\left|\sqrt{1+e_2^{-4}\tan^2\theta+e_1^{-4}e_2^{-4}+e_1^{-4}\tan^2\theta}-1\right|_{s,(1+\eta_1)K}\leq \lambda^{-\frac{1}{12}q_{n-1}}.
\end{align}
\begin{align}\label{ne3}
\left|\sqrt{\frac{1}{1+e_2^{-4}\tan^2\theta+e_1^{-4}e_2^{-4}+e_1^{-4}\tan^2\theta}}-1\right|_{s,(1+\eta_2)(1+\eta_1)K}\leq \lambda^{-\frac{1}{144}q_{n-1}}.
\end{align}
By (1) in Proposition \ref{progev}, \eqref{lowcond2}, \eqref{lowcond1}, \eqref{ne2}, \eqref{ne9} and \eqref{ne3}, we have
\begin{align}\label{ne8}
\nonumber \left|\frac{1}{b}\right|_{s,(1+\eta_1)K}&\leq |e_1^{-1}|_{s,K}^{2}|e_2^{-1}|_{s,K}^2\left|\frac{1}{\cos^2\theta}\right|_{s,K}\left|\frac{1}{1+e_2^{-4}\tan^2\theta+e_1^{-4}e_2^{-4}+e_1^{-4}\tan^2\theta}\right|_{s,(1+\eta_1)K}\\
&\leq C^2e^{2q_n^\gamma}\lambda^{-\frac{4}{3}q_{n-1}}(1+\lambda^{-\frac{1}{12}q_{n-1}})\leq \lambda^{-q_{n-1}}.
\end{align}
\begin{align}\label{ne7}
\nonumber \left|\sqrt{b}\right|_{s,(1+\eta_1)K}&\leq \left|e_1\right|_{s,K}\left|e_2\right|_{s,K}\left|\cos\theta\right|_{s,K}\left|\sqrt{1+e_2^{-4}\tan^2\theta+e_1^{-4}e_2^{-4}+e_1^{-4}\tan^2\theta}\right|_{s,(1+\eta_1)K}\\
&\leq 2C\left|e_1\right|_{s,K}\left|e_2\right|_{s,K}.
\end{align}
\begin{align}\label{ne4}
\nonumber \left|\sqrt{\frac{1}{b}}\right|_{s,(1+\eta_2)(1+\eta_1)K}&\leq \left|\frac{1}{e_1}\right|_{s,K}\left|\frac{1}{e_2}\right|_{s,K}\left|\frac{1}{\cos\theta}\right|_{s,K}\left|\sqrt{\frac{1}{1+e_2^{-4}\tan^2\theta+e_1^{-4}e_2^{-4}+e_1^{-4}\tan^2\theta}}\right|_{s,(1+\eta_2)(1+\eta_1)K}\\
&\leq 2C|e^{-1}_1|_{s,K}|e^{-1}_2|_{s,K}e^{q_n^\gamma}.
\end{align}
By \eqref{ne8}, (3) and (4) in Proposition \ref{progev}, we have
$$
\left|\sqrt{1-4b^{-2}}-1\right|_{s,(1+\eta_1)^2K}\leq \lambda^{-\frac{1}{12}q_{n-1}},
$$
\begin{align}\label{ne10}
\left|\sqrt{\frac{1+\sqrt{1-4b^{-2}}}{2}}-1\right|_{s,(1+\eta_2)(1+\eta_1)^2K}\leq  \lambda^{-\frac{1}{12^2}q_{n-1}},
\end{align}
\begin{align}\label{ne5}
\left|\sqrt{\frac{2}{1+\sqrt{1-4b^{-2}}}}-1\right|_{s,(1+\eta_3)(1+\eta_2)(1+\eta_1)^2K}\leq  \lambda^{-\frac{1}{12^3}q_{n-1}},
\end{align}
where $\eta_2=\lambda^{-\frac{1}{2880}q_{n-1}}$.

By \eqref{ne7}, \eqref{ne4}, \eqref{ne10} and \eqref{ne5}, we have
\begin{align*}
|e|_{s,(1+\eta)K}=\left|\sqrt{{b}}\right|_{s,(1+\eta)K}\left|\sqrt{\frac{1+\sqrt{1-4b^{-2}}}{2}}\right|_{s,(1+\eta)K}\leq C^2\left|e_1\right|_{s,K}\left|e_2\right|_{s,K}.
\end{align*}
\begin{align*}
|e^{-1}|_{s,(1+\eta)K}=\left|\sqrt{\frac{1}{b}}\right|_{s,(1+\eta)K}\left|\sqrt{\frac{2}{1+\sqrt{1-4b^{-2}}}}\right|_{s,(1+\eta)K}\leq C^2|e^{-1}_1|_{s,K}|e^{-1}_2|_{s,K}e^{q_n^\gamma}.
\end{align*}
We finish the proofs of \eqref{r3}-\eqref{r1}.

Now we prove (\ref{r2}).  By polar decomposition procedure, we have  $s(x)=\frac{\pi}{2}+\theta(E(x))$ where $s(x)$ is the most contraction direction of $E(x)$ and $\theta(E(x))$ is the eigen-direction of $E^t(x)E(x)$ corresponding to the eigenvalue $\|E(x)\|^2$.  Let $a=e_1e_2$, $c=\frac{e_1}{e_2}$, $u=2(e_2^2-e_2^{-2})\sin\theta\cos\theta$ and $U=(a^2-a^{-2})\cos^2\theta+(c^2-c^{-2})\sin^2\theta$. It's easy to calculate that
\begin{equation}\label{a1}
\tan s(x)=\tan (\frac{\pi}{2}+\theta(E(x)))=\frac{u}{\sqrt{U^2+u^2}-U}=\frac{\sqrt{U^2+u^2}+U}{u}.
\end{equation}
Since $|\theta-\pi/2|>ce^{-q_n^\gamma}\gg  \min\left\{\inf\limits_{x\in I} e_1(x),\inf\limits_{x\in I} e_2(x)\right\}^{-\frac{1}{100}}$, we have
\begin{equation}\label{a2}
U\geq ce_1^2e_2^2e^{-2q_n^\gamma}-e_1^2e_2^{-2}-e_1^{-2}e_2^{2}>0.
\end{equation}
\eqref{a1} and \eqref{a2} imply that
\begin{equation}\label{a3}
s(x)=\arctan\left(sgn(u)\left(\sqrt{\frac{U^2}{u^2}+1}+\frac{U}{|u|}\right)\right).
\end{equation}
A direct calculation shows that
\begin{align}\label{a4}
g:=\frac{U}{u}=\frac{(e_1e_2)^2-(e_1e_2)^{-2}}{2(e_2^2-e_2^{-2})}\cot\theta+\frac{(e_1e^{-1}_2)^2-(e_1e^{-1}_2)^{-2}}{2(e_2^2-e_2^{-2})}\tan\theta.
\end{align}
Without loss of generality, we only consider the case $g(x)>0$. A direct computation shows
\begin{align*}
\frac{ds}{dx}=\frac{1}{2}\frac{1}{1+g^2}\frac{dg}{dx},
\end{align*}
$$
g=\frac{e_1^2\cot\theta}{2}\frac{1-e_1^{-2}e_2^{-4}+(e_2^{-4}-e_1^{-4})\tan^2\theta}{1-e_2^{-4}}:=\frac{e_1^2\cot\theta}{2}\cdot h,
$$
\begin{align*}
\frac{1}{1+g^2}=4e_1^{-4}\tan^2\theta\frac{1}{4e_1^{-4}\tan^2\theta+h^2}.
\end{align*}
By \eqref{lowcond2}, \eqref{lowcond1}, (1) and (3) in Proposition \ref{progev}, we have
$$
|h-1|_{s,(1+\eta_1)K}=\left|\frac{1-e_1^{-2}e_2^{-4}+(e_2^{-4}-e_1^{-4})\tan^2\theta}{1-e_2^{-4}}-1\right|_{s,(1+\eta_1)K}\leq \lambda^{-q_{n-1}},
$$
$$
|h^2-1|_{s,(1+\eta_1)K}\leq |h-1|_{s,(1+\eta_1)K}|h+1|_{s,(1+\eta_1)K}\leq 4\lambda^{-q_{n-1}},
$$
$$
\left|\frac{1}{4e_1^{-4}\tan^2\theta+h^2}-1\right|_{s,(1+\eta_1)^2K}\leq 4\lambda^{-\frac{q_{n-1}}{12}}.
$$
It follows that
$$
|g|_{s,(1+\eta_1)K}\leq C|e_1|^2_{s,K}|\cot\theta|_{s,K}\left|\frac{1-e_1^{-2}e_2^{-4}+(e_2^{-4}-e_1^{-4})\tan^2\theta}{1-e_2^{-4}}\right|_{s,(1+\eta_1)K}\leq C|e_1|^2_{s,K},
$$
\begin{align*}
\left|\frac{1}{1+g^2}\right|_{s,(1+\eta_1)^2K}\leq 2|\tan\theta|^2_{s,K}|e_1^{-1}|^{4}_{s,K}\leq C|e_1^{-1}|^{4}_{s,K}e^{2q_n^\gamma}.
\end{align*}
By (2) in Proposition \ref{progev} and \eqref{lowcond2}, we have
\begin{equation*}
\left|\frac{\pi}{2}-s\right|_{s,(1+\eta)K}\leq \left|\frac{1}{1+g^2}\right|_{s,(1+\eta_1)^2K}|\partial g|_{s,(1+\eta)K}\leq  C^2|e_1^{-1}|^{4}_{s,K}e^{2q_n^\gamma}|e_1|^2_{s,K}K\lambda^{\frac{1}{20}q_{n-1}}\leq |e_1^{-1}|^{3}_{s,K}|e_1|^2_{s,K}.
\end{equation*}
Similar results hold for $u$. Thus we finish the whole proof.
\end{pf}
Consider a sequence of maps
$$
E^\ell\in G^{s,K}(I,SL(2,\R)), \ \ 0\leq \ell\leq n-1.
$$
Let $s^\ell(x)=s[E^\ell(x)]$, $u^\ell(x)=u[E^\ell(x)]$, $e^\ell(x)=\|E^\ell(x)\|$ and $\Lambda^\ell(x)=\begin{pmatrix}e^\ell(x)&0\\ 0&(e^\ell(x))^{-1}\end{pmatrix}$. By polar decomposition, it holds that
$$
E^\ell(x)=R_{u^\ell(x)}\Lambda^\ell(x)R_{\frac{\pi}{2}-s^\ell(x)}.
$$
Set for each $0\leq \ell\leq n-1$,
\begin{align*}
\begin{split}
E_k^\ell(x)=\left\{
\begin{array}{ll}
E^{k-1+\ell}(x)\cdots E^\ell(x), & 1\leq k\leq n-\ell\\
I_2, & k=0\\
\left(E_{-k}^{\ell+k}(x)\right)^{-1}, & -\ell\leq k\leq -1
\end{array}
\right.
\end{split}.
\end{align*}
For $k\geq 1$, let $s_k^\ell(x)=s[E_k^\ell(x)]$, $u_k^\ell(x)=s[E_{-k}^\ell(x)]$, $e_k^\ell(x)=\|E_k^\ell(x)\|$ and $\Lambda_k^\ell(x)=\begin{pmatrix}e_k^\ell(x)&0\\ 0&(e_k^\ell(x))^{-1}\end{pmatrix}$. Again by polar decomposition, it holds that
$$
E_k^\ell(x)=R_{u_k^{k+\ell}(x)}\Lambda_k^\ell(x)R_{\frac{\pi}{2}-s_k^\ell(x)}.
$$
\begin{Lemma}\label{2}
Let  $0\leq \ell<n-1\leq q_m^C$, $1<\xi<\frac{1}{1000}$ and $\theta_\ell(x)=:u^{\ell-1}(x)-s^\ell(x)+\frac{\pi}{2}$. Assume that
\begin{equation}\label{nlowcond}
\inf\limits_{x\in I}\left|\frac{\pi}{2}-\theta_\ell(x)\right|=\inf\limits_{x\in I}\left|s^\ell(x)-u^{\ell-1}(x)\right|\geq ce^{-q_m^\gamma}\gg\min_{0\leq \ell\leq n-1}\left\{\inf_{x\in I}e^\ell(x)\right\}^{-\frac{1}{100}},
\end{equation}
\begin{equation}\label{nlowcond1}
 \ \ \left|\frac{1}{\cos \theta_\ell}\right|_{s,K}, \left|\tan \theta_\ell\right|_{s,K}\leq Ce^{q_m^\gamma}, \ \ \left|\cos\theta_\ell\right|_{s,K},\left|\cot\theta_\ell\right|_{s,K}\leq C,
\end{equation}
\begin{equation}\label{nlowcond2}
|(e^\ell)^{-1}|_{s,K}\leq \left(|e^\ell|_{s,K}\right)^{-1+\xi}\leq C\lambda^{-\frac{1}{3}q_{m-1}}.
\end{equation}
Then it holds that
\begin{equation}\label{nr3}
\inf\limits_{x\in I}e_n^0(x)\geq c^ne^{-2nq_m^\gamma}\prod\limits_{\ell=0}^{n-1}\inf\limits_{x\in I}e^\ell(x),
\end{equation}
\begin{equation}\label{nr4}
|e_n^0|_{s,(1+\eta)K}\leq C^{2n}\prod\limits_{\ell=0}^{n-1}|e^\ell|_{s,K},
\end{equation}
\begin{equation}\label{nr1}
|(e_n^0)^{-1}|_{s,(1+\eta)K}\leq C^{2n}\prod\limits_{\ell=0}^{n-1}|(e^\ell)^{-1}|_{s,K}e^{4nq_m^\gamma},
\end{equation}
\begin{equation}\label{nr2}
\left|s^0-s_n^0\right|_{s,(1+\eta)}\leq \lambda^{-\frac{1}{10}q_{m-1}},\ \ \left|u^{n-1}-u_n^n\right|_{s,(1+\eta)K}\leq \lambda^{-\frac{1}{10}q_{m-1}},
\end{equation}
where $\eta=\lambda^{-\frac{1}{8000}q_{m-1}}$.
\end{Lemma}
\begin{pf}
We prove it by induction. In the case of the product of two matrices, i.e. $n=2$, it follows from Lemma \ref{2}. Now, we assume for $k\leq n-1$ and for all possible $\ell$, we have
\begin{equation}\label{inr3}
\inf\limits_{x\in I}e_k^\ell(x)\geq c^ke^{-2kq_m^\gamma}\prod\limits_{j=0}^{k-1}\inf\limits_{x\in I}e^{\ell+j}(x),
\end{equation}
\begin{equation}\label{inr4}
|e_k^\ell|_{s,(1+\eta)^{2k}K}\leq C^{2k}\prod\limits_{j=0}^{k-1}|e^{\ell+j}|_{s,K},
\end{equation}
\begin{equation}\label{inr1}
|(e_k^\ell)^{-1}|_{s,(1+\eta)^{2k}K}\leq C^{2k}\prod\limits_{j=0}^{k-1}|(e^{\ell+j})^{-1}|_{s,K}e^{4kq_m^\gamma},
\end{equation}
\begin{equation}\label{inr2}
\left|s_{k-1}^\ell-s_k^\ell\right|_{s,(1+\eta)^{2k}K}\leq \lambda^{-\frac{1}{9}kq_{m-1}},
\end{equation}
\begin{equation}
\left|u_{k-1}^{n-\ell}-u_k^{n-\ell}\right|_{s,(1+\eta)^{2k}K}\leq \lambda^{-\frac{1}{9}kq_{m-1}},
\end{equation}
where $\eta=\lambda^{-\frac{1}{4000}q_{m-1}}$.

Clearly, \eqref{inr2} implies that for $\ell=0,1$,
\begin{align*}
\left|s^\ell-s_{n-1}^\ell\right|_{s,(1+\eta)^{2n-2}K}\leq \sum_{j=1}^{n-1}\lambda^{-\frac{1}{9}jq_{m-1}}\leq 2\lambda^{-\frac{1}{9}q_{m-1}},
\end{align*}
\begin{align*}
\left|u^{n-\ell-1}-u_{n-1}^{n-\ell}\right|_{s,(1+\eta)^{2n-2}K}\leq \sum_{j=1}^{n-1}\lambda^{-\frac{1}{9}jq_{m-1}}\leq 2\lambda^{-\frac{1}{9}q_{m-1}}.
\end{align*}
Combining the above with \eqref{nlowcond}, we have
\begin{equation}\label{need1}
\inf\limits_{x\in I}\left|u_{n-1}^{n-1}(x)-s^{n-1}(x)\right|\geq ce^{-q_m^\gamma},
\end{equation}
\begin{equation*}
\inf\limits_{x\in I}\left|s_{n-1}^1(x)-u^0(x)\right|\geq ce^{-q_m^\gamma}.
\end{equation*}
 Let $\tilde{\theta}_{n-1}(x)=u_{n-1}^{n-1}(x)-s^{n-1}(x)$, then
\begin{align*}
\left|\tilde{\theta}_{n-1}-\theta_{n-1}\right|_{G^{s,(1+\eta)^{2n-2}K}(I)}\leq \left|u_{n-1}^{n-1}-u^n\right|_{G^{s,(1+\eta)^{2n-2}K}(I)} \leq 2\lambda^{-\frac{1}{9}q_{m-1}}.
\end{align*}
Note that
\begin{align*}
\frac{1}{\cos \tilde{\theta}_{n-1}}=\frac{1}{\cos\theta_{n-1}}\frac{1}{\cos(\tilde{\theta}_{n-1}-\theta_{n-1})-\tan\theta_{n-1}\sin(\tilde{\theta}_{n-1}-\theta_{n-1})}.
\end{align*}
By (5) in Proposition \ref{progev} and \eqref{nlowcond1}, we have
$$
\left|{\cos(\tilde{\theta}_{n-1}-\theta_{n-1})-\tan\theta_{n-1}\sin(\tilde{\theta}_{n-1}-\theta_{n-1})}-1\right|_{G^{s,(1+\eta)^{2n-1}K}}\leq  \lambda^{-\frac{1}{120}q_{m-1}}.
$$
By (3) in Proposition \ref{progev}, we have
\begin{equation}\label{need2}
\left|\frac{1}{\cos(\tilde{\theta}_{n-1}-\theta_{n-1})-\tan\theta_{n-1}\sin(\tilde{\theta}_{n-1}-\theta_{n-1})}\right|_{G^{s,(1+\eta)^{2n-1}K}}\leq 2.
\end{equation}
By \eqref{nlowcond1} and \eqref{need2}, we have
\begin{align}\label{need3}
 \left|\frac{1}{\cos \tilde{\theta}_{n-1}}\right|_{G^{s,(1+\eta)^{2n-1}K}(I)}\leq 2 \left|\frac{1}{\cos \theta_{n-1}}\right|_{G^{s,(1+\eta)^{2n-1}K}(I)}\leq 2Ce^{q_m^\gamma}.
\end{align}
Similarly
\begin{align}\label{need4}
\left|\tan \tilde{\theta}_{n-1}\right|_{G^{s,(1+\eta)^{2n-1}K}(I)}\leq 2Ce^{q_m^\gamma}, \left|\cos\tilde{\theta}_{n-1}\right|_{G^{s,(1+\eta)^{2n-1}K}(I)},\left|\cot\tilde{\theta}_{n-1}\right|_{G^{s,(1+\eta)^{2n-1}K}(I)}\leq 2C.
\end{align}
\eqref{inr4}, \eqref{inr1}, \eqref{need1}, \eqref{need3} and \eqref{need4} imply that we can apply Lemma \ref{2} to the product
$$
E_n^0(x)=E^{n-1}(x)E_{n-1}^0(x)=E_{n-1}^1(x)E^0(x),
$$
which implies that
\begin{equation*}
\inf\limits_{x\in I}e_n^0(x)\geq c\inf\limits_{x\in I}e^0_{n-1}(x)\inf\limits_{x\in I}e^{n-1}(x)e^{-q_m^\gamma}\geq c^ne^{-2nq_m^\gamma}\prod\limits_{\ell=0}^{n-1}\inf\limits_{x\in I}e^{\ell}(x),
\end{equation*}
\begin{equation*}
|e_n^0|_{s,(1+\eta)^{2n}K}\leq C^2|e^{n-1}|_{s,(1+\eta)^{2(n-1)}K}|e_{n-1}^0|_{s,(1+\eta)^{2(n-1)}K}\leq C^{2n}\prod\limits_{\ell=0}^{n-1}|e^\ell|_{s,K},
\end{equation*}
\begin{equation*}
|(e_n^0)^{-1}|_{s,(1+\eta)^{2n}K}\leq C^2e^{2q_m^\gamma}|(e^{n-1})^{-1}|_{s,(1+\eta)^{2(n-1)}K}|(e_{n-1}^0)^{-1}|_{s,(1+\eta)^{2(n-1)}K}\leq C^{2n}e^{4nq_m^\gamma}\prod\limits_{\ell=0}^{n-1}|(e^\ell)^{-1}|_{s,K}.
\end{equation*}
By Lemma  \ref{2} and \eqref{nlowcond2}, we have
\begin{align*}
\left|s_{n-1}^0-s_n^0\right|_{s,(1+\eta)^{2n}K}&\leq |(e_{n-1}^0)^{-1}|^{3}_{s,(1+\eta)^{2n-2}K}|e_{n-1}^0|^2_{s,(1+\eta)^{2n-2}K}\\
&\leq C^{6n-6}e^{12nq_m^\gamma}\prod\limits_{\ell=0}^{n-2}|(e^\ell)^{-1}|^3_{s,(1+\eta)^{2n-2}K} C^{4n-4}\prod\limits_{\ell=0}^{n-2}|e^\ell|^2_{s,(1+\eta)^{2n-2}K}\\
&\leq C^{10n}e^{12nq_m^\gamma}\prod\limits_{j=0}^{n-2}|(e^\ell)^{-1}|^\frac{1}{2}_{s,(1+\eta)^{2n-2}K}\leq \lambda^{-\frac{1}{9}nq_{m-1}},
\end{align*}
thus
\begin{equation*}
\left|s^0-s_n^0\right|_{s,(1+\eta)^{2n}K}\leq \left|s^0-s_{n-1}^0\right|_{s,(1+\eta)^{2n}K}+\left|s_{n-1}^0-s_n^0\right|_{s,(1+\eta)^{2n}K}\leq \lambda^{-\frac{1}{10}q_{m-1}}.
\end{equation*}
Similar results hold for $u_n^n$, we finish the proof since $(1+\eta)^{2n}\leq 1+\lambda^{-\frac{1}{8000}q_{m-1}}$ for $n<q_m^C$.
\end{pf}

In the following, we will fix $0<\nu_2<\nu_1<1$, $s_2=1+\frac{1}{\nu_2}>s_1=1+\frac{1}{\nu_1}>2$, $\beta>1$ such that $0<\gamma_2=\beta\nu_2<\gamma_1=\beta\nu_1<1$. Let $\delta_1>0$ be sufficiently small such that $0<\frac{\nu_1\beta}{1-\delta_1\nu_1}<1$.  Recall that
\begin{itemize}
\item The critical set: $C_0=\{c_1,c_2\}$ where $c_1\in [0,\pi)$ and $c_2=c_1+\pi$.
\item The critical interval: $I_{n,1}=\left[c_1-\frac{1}{q_n^\beta},c_1+\frac{1}{q_n^\beta}\right]$, $I_{n,2}=\left[c_2-\frac{1}{q_n^\beta},c_2+\frac{1}{q_n^\beta}\right]$ and $I_n=I_{n,1}\bigcup I_{n,2}$.
\item The first return time: For $x\in I_n$, we denote the smallest positive integer $i$ with $T^ix\in I_n$ (respectively $T^{-i}x\in I_n$) by $r_n^+(x)$ (respectively $r_n^-(x)$), and define $r_n^{\pm}=\min_{x\in I_n}r_n^{\pm}(x)$. Obviously, $r_n^{\pm}\geq\frac{q_n}{2}$.
\end{itemize}
\begin{Remark}
If $\alpha$ is bounded, we have $r_n^{\pm}\leq q_n^C$ for some $C$ only depending on $\alpha$. See \cite{jk} for the proof.
\end{Remark}

\subsection{Proof of Proposition \ref{p1}}We prove Proposition \ref{p1} by induction. Instead of $r_n^{\pm}(x)$, sometimes, we use $r_n^{\pm}$ for short when the difference between $r_n^{\pm}(x)$ and $r_n^{\pm}$ are negligible.  Recall that
$$
\ln\lambda_{n+1}=\ln\lambda_n-10^4 q_{n+1}^{\gamma_1-1},\ \ \gamma_1=\nu_1\beta, \ \ \lambda_N=\lambda^{1-\e}.
$$
$$
\ln\widetilde{\lambda_{n+1}}=\ln\widetilde{\lambda_n}+10^4 q_{n+1}^{\gamma_1-1},\ \ \gamma_1=\nu_1\beta, \ \ \widetilde{\lambda_N}=\lambda^{1+\e}.
$$

We first construct $\phi_N(x)$ and $A_N(x)$ such that $(1)_N - (4)_N$ hold.

\noindent {\it Construction of $\phi_N(x)$ and $A_N(x)$}: Let $c_1,c_2\in \T$ with $c_1\in [0,\pi)$ and $c_2=c_1+\pi$. We define a $2\pi$-periodic smooth function $\phi_0$ by
$$
\sin(\phi_0(x))=ce^{-\left(\frac{1}{(x-c_1-k\pi)^{\nu_1}}+\frac{1}{(c_1+(k+1)\pi-x)^{\nu_1}}\right)},\ \ x\in [c_1+k\pi,c_1+(k+1)\pi),
$$
for some $0<c<\frac{1}{1000}$. In view of Proposition \ref{progev} and Corollary \ref{phi_0}, it's easy to see
\begin{enumerate}
\item $\phi_0$ is a $G^{1+\frac{1}{\nu_1},C}$-$2\pi$ periodic function for some $C>0$.
\item $|\phi_0|_{C^0(\mathbb{S}^1)}\leq \frac{\pi}{6}$ and for $i=1,2$, $|\phi_0(x)|\geq ce^{-{|x-c_i|^{-\nu_1}}}$ for some $c>0$.
\end{enumerate}
Let $A(x)=\Lambda\cdot R_{\frac{\pi}{2}-\phi_0(x)}=\begin{pmatrix}\lambda&0\\ 0&\lambda^{-1}\end{pmatrix}\cdot R_{\frac{\pi}{2}-\phi_0(x)}$, by \cite{young}, there exists a large $\lambda_0>0$ depending on $\phi_0$, $\nu_1$ and $\e$ such that if $\lambda>\lambda_0$,
$$
\{A(x),\cdots,A(T^{r_N^+(x)-1}x)\}\ \ \text{is}\ \ \lambda_N-\text{hyperbolic}, \forall x\in I_N.
$$

Let $\overline{s}_N(x)=\overline{s(A^{r_N^+}(x))}$, $\overline{s}_N^\prime(x)=\overline{s(A^{-r_N^-}(x))}$ for $x\in I_N$. Let $e_N(x)$ be a $2\pi$-periodic $C^\infty$-function such that $e_N(x)=\phi_0(x)-(\overline{s}_N^\prime(x)-\overline{s}_N(x))$ for $x\in I_N$. Let $\hat{e}_N(x)=e_N(x)\cdot f_N(x)$ where $f_N$ is defined in Lemma \ref{fn} and $\phi_N(x)=\phi(x)+\hat{e}_N(x)$ for $x\in \mathbb{S}^1$. \\
\noindent
{\it  Verifying  $(1)_N$ and $(4)_N$ of Proposition \ref{p1}}: Let $\eta_N=\lambda_N^{-\frac{1}{8000}}$, $n=r_N^+\leq q_N^C$, $e^\ell(x)=\|A(x+\ell\alpha)\|$ and $I=I_N$,
$$
E^\ell(x)=A(x+\ell\alpha)=\Lambda\cdot R_{\frac{\pi}{2}-\phi_0(x+\ell\alpha)}.
$$
For $0\leq \ell<n-1$, since $x+\ell\alpha\notin I_N$, one can easily verify that
$$
\inf_{x\in I}\|A(x+\ell\alpha)\|=\lambda,
$$
\begin{equation}\label{nd1}
\inf\limits_{x\in I}\left|\frac{\pi}{2}-\theta_\ell(x)\right|:=\inf\limits_{x\in I}\left|\phi_0(x+\ell\alpha)\right|\geq ce^{-q_N^{\beta\nu_1}}=ce^{-q_N^{\gamma_1}},
\end{equation}
\begin{equation}\label{nd2}
|(e^\ell)^{-1}|_{G^{s_1,C}(I)}=|e^\ell|_{G^{s_1,C}(I)}^{-1}=\lambda^{-1}.
\end{equation}
By Corollary \ref{inverse-Gevrey-function} and Corollary \ref{phi_0}, there is some $C>0$ such that
\begin{equation}\label{nd4}
|\cos\theta_\ell|_{G^{s_1,C}(I)}=\left|\sin(\phi_0(x+\ell\alpha))\right|_{G^{s_1,C}(I)},\ \  \left|\cot\theta_\ell\right|_{G^{s_1,C}(I)}=\left|\tan(\phi_0(x+\ell\alpha))\right|_{G^{s_1,C}(I)}\leq C,
\end{equation}
\begin{equation}\label{nd5}
\left|\frac{1}{\cos\theta_\ell}\right|_{G^{s_1,C}(I)}=\left|\frac{1}{\sin(\phi_0(x+\ell\alpha))}\right|_{G^{s_1,C}(I)},\ \  \left|\tan\theta_\ell\right|_{G^{s_1,C}(I)}=\left|\cot(\phi_0(x+\ell\alpha))\right|_{G^{s_1,C}(I)}\leq Ce^{q_N^{\gamma_1}}.
\end{equation}
Set $q_{N-1}=1$, \eqref{nd1}-\eqref{nd5} imply that all the assumptions in Lemma \ref{2} are satisfied. It follows
$$
\inf_{x\in I_N}\|A_{r_N^+}(x)\|\geq \lambda^{r_N^+}c^{r_N^+}e^{-r_N^+q_N^{\gamma_1}}\geq \lambda_N^{r_N^+}, \ \  \|A_{r_N^+}\|_{G^{s_1,(1+\eta_N)C}(I_N)}\leq C^{2r_N^+}\lambda^{r_N^+}\leq \widetilde{\lambda_N}^{r_N^+},
$$
$$
\left|\frac{1}{\|A_{r_N^+}\|}\right|_{G^{s_1,(1+\eta_N)C}(I_N)}\leq C^{2r_N^+}\lambda^{-r_N^+}e^{4r_N^+q^{\gamma_1}_{N}}\leq \lambda_N^{-r_N^+}, \ \ |e_N|_{G^{s_1,(1+\eta_N)C}(I)}\leq \lambda_N^{-\frac{1}{10}}.
$$
Similar results hold for $r_N^-$, we omit the proof.

By the definition and Lemma \ref{fn}, we have
\begin{align*}
|\phi_N-\phi|_{G^{s_1,(1+\eta_N)C}(\mathbb{S}^1)}&=|\hat{e}_N|_{G^{s_1 ,(1+\eta_N)C}(\mathbb{S}^1)}\leq |e_N|_{G^{s_1,(1+\eta_N)C}(I_N)}|f_N|_{G^{s_1,(1+\eta_N)C}(\mathbb{S}^1)}\\
&\leq \lambda_N^{-\frac{1}{10}}(Cq^\beta_N)^{q_N^{\frac{\nu_1\beta}{1-\delta_1\nu_1}}}\leq \lambda_N^{-\frac{1}{20}}.
\end{align*}
The last inequality holds since $\frac{\nu_1\beta}{1-\delta_1\nu_1}<1$ and $\lambda\gg e^{q_N^{q_N}}$.\\

\noindent { \it Verifying  $(2)_N$ of  Proposition \ref{p1}}: Let $A_N(x)=\Lambda\cdot R_{\frac{\pi}{2}-\phi_N(x)}$. Obviously, $A_N(x)=A(x)\cdot R_{-\hat{e}_N(x)}$.
\begin{Lemma}[\cite{wangyou}]
For $x\in I_N$, it holds that
$$
A_N^{r_N^+}(x)=A^{r_N^+}(x)\cdot R_{-\hat{e}_N(x)}
$$
and
$$
A_N^{-r_N^-}(x)=R_{\hat{e}_N(T^{-r_N^-}x)}\cdot A^{-r_N^-}(x).
$$
\end{Lemma}
Thus, for any $x\in I_N$, $\{A_N(x),...,A_N(T^{r_N^+(x)-1}x)\}$ is a $\lambda_N$-hyperbolic sequence.\\

\noindent
{\it Verifying  $(3)_N$ of Proposition \ref{p1}}:  $(s_N-s_N^\prime)(x)=(\overline{s}_N-\overline{s}_N^\prime)(x)+\hat{e}_N(x)$ which implies that $(s_N-s_N^\prime)(x)=\phi_0(x)$ on $\frac{I_N}{10}$, since $|e_N(x)|\leq \lambda_N^{-\frac{1}{20}}$ in $I_N$. Thus we have
$$
|(s_N-s_N^\prime)(x)|\geq |\phi_0(x)|-\lambda_N^{-\frac{1}{20}}\geq ce^{-10^{\nu_1} q_N^{\gamma_1}},
$$
on $I_N\backslash\frac{I_N}{10}$ since $\lambda>e^{q_N^{q_N}}$.\\

Inductively, we assume that $\phi_N(x),...,\phi_{n-1}(x)$ have been constructed such that Proposition \ref{p1} holds for $N\leq i\leq n-1$, i.e., \\
\begin{enumerate}[$(1)_{i}$]
\item $ |\phi_i(x)-\phi_{i-1}(x)|_{s_1,(1+\eta_i)C}\leq \lambda_i^{-\frac{q_{i-1}}{100}}$ where $\eta_i=\prod_{j=N}^{i-1}(1+\lambda_j^{-\frac{1}{8000}q_{j}})-1$. \\

\item For each $x\in I_i$, $A_i(x),A_i(Tx),\cdots,A_i(T^{r_i^+(x)-1}(x))$ is $\lambda_i$-hyperbolic.\\

\item We have
\begin{align*}
(a)_i\ \ s_i(x)-s_i^\prime(x)=\phi_0(x)\ \ x\in \frac{I_i}{10};
\end{align*}
\begin{align*}
(b)_i\ \ |s_i(x)-s_i^\prime(x)|\geq \frac{1}{2}|\phi_0(x)|\geq \frac{1}{2}e^{-10^{\nu_1} q_i^{\nu_1\beta}}, \ \ x\in I_i\backslash \frac{I_i}{10}.
\end{align*}

\item  It holds
$$
\|A_{r_i^{\pm}}\|_{G^{s_1,(1+\eta_i)C}(I_n)}\leq \widetilde{\lambda_i}^{r_i^{\pm}},\ \ \left|\frac{1}{\|A_{r_i^{\pm}}\|}\right|_{G^{s_1,(1+\eta_i)C}(I_n)}\leq \lambda_i^{-r_i^{\pm}}.
$$\\
\end{enumerate}
Now we construct $\phi_n(x)$ and verify $(1)_n - (4)_n$. \\

\noindent {\it Constructing $\phi_n(x)$}: From $(2)_{n-1}$, we have that
$$
\left\|A_{n-1}^{r_{n-1}^+(x)}(x)\right\|\cdot e^{-(10q_{n-1}^\beta)^{\nu_1}}\geq \lambda_{n-1}^{q_{n-1}}\cdot e^{-(10q_{n-1}^\beta)^{\nu_1}}\geq \lambda_n^{(1-\epsilon)q_{n-1}},\ \ x\in I_{n-1}.
$$
Combing the above with $(3)_{n-1}$, for each $x\in I_{n}$, $A_{n-1}(x),A_{n-1}(Tx),\cdots,A_{n-1}(T^{r_n^+(x)-1}(x))$ is $\lambda_{n}$-hyperbolic.
Let $\overline{s}_n(x)=\overline{s(A^{r_n^+}(x))}$, $\overline{s}_n^\prime(x)=\overline{s(A^{-r_n^-}(x))}$. Let $e_n(x)$ be a $2\pi$-periodic $C^\infty$-function such that $e_n(x)=\phi_0(x)-(\overline{s}_n^\prime(x)-\overline{s}_n(x))$ for $x\in I_n$. Let $\hat{e}_n(x)=e_n(x)\cdot f_n(x)$ and $\phi_n(x)=\phi(x)+\hat{e}_n(x)$ for $x\in \mathbb{S}^1$. \\

\noindent {\it Verifying  $(1)_n$ and $(4)_n$ of Proposition \ref{p1}}: For any $x\in I_n$,  let $j_i$ be defined so that $T^{j_i}x\in I_{n-1}\backslash I_{n}$ and let $T^{j_{i+1}}x$ be the next return of $T^{j_i}x$ to $I_{n-1}$. Let $\eta_n=\prod_{j=N}^{n-1}(1+\lambda_j^{-\frac{1}{8000}q_{j}})-1$, $n_0\leq q_{n}^C$, $e^\ell(x)=\|A_{j_{\ell+1}-j_\ell}(x+j_\ell\alpha)\|$ and $I=I_n$.
$$
E^\ell(x)=A_{j_{\ell+1}-j_\ell}(x+j_\ell\alpha)=R_{u_{j_{\ell+1}-j_\ell}(x+j_{\ell+1}\alpha)}\begin{pmatrix}e_\ell(x)&0\\ 0&(e_\ell(x))^{-1}\end{pmatrix}R_{\frac{\pi}{2}-s_{j_{\ell+1}-j_\ell}(x+j_\ell\alpha)}.
$$
By $(2)_{n-1}$, we have
$$
\inf_{x\in I}\|A_{j_{\ell+1}-j_\ell}(x+j_\ell\alpha)\|\geq \lambda_{n-1}^{j_{\ell+1}-j_\ell}\geq \lambda_{n-1}^{\frac{q_{n-1}}{2}}, \ \ 0\leq \ell\leq n_0-1.
$$
By $(3)_{n-1}$, for $0\leq \ell<n_0-1$,
\begin{align*}
\inf\limits_{x\in I}\left|\frac{\pi}{2}-\theta_\ell(x)\right|&:=\left|s_{j_{\ell+1}-j_\ell}(x+j_\ell\alpha)-u_{j_{\ell}-j_{\ell-1}}(x+j_\ell\alpha)\right|\\
&\geq \frac{1}{2}\left|s_{n-1}(x+j_\ell\alpha)-s'_{n-1}(x+j_\ell\alpha)\right|\geq ce^{-q_{n}^{\beta\nu_1}}=ce^{-q_{n}^{\gamma_1}}.
\end{align*}
By $(4)_{n-1}$, we have
$$
|e^\ell|_{G^{s,(1+\eta_{n-1})C}}\left|\frac{1}{e^\ell}\right|_{G^{s,(1+\eta_{n-1})C}}\leq\left(\frac{\widetilde{\lambda_{n-1}}}{\lambda_{n-1}}\right)^{j_{\ell+1}-j_\ell}\leq \lambda^{4\e({j_{\ell+1}-j_\ell})}\leq |e_\ell|_{G^{s,(1+\eta_{n-1})C}}^{\xi},
$$
$$
|e^\ell|_{G^{s,(1+\eta_{n-1})C}}^{-1+\xi}\leq \lambda_{n-1}^{-\frac{j_{\ell+1}-j_\ell}{2}}\leq  \lambda_{n-1}^{-\frac{q_{n-1}}{3}}.
$$
By $(1)_{n-1}$, we have $|\phi_{n-1}-\phi_0|_{G^{s,(1+\eta_{n-1})C}}\leq 2\lambda^{-\frac{1}{100}}$. By Proposition \ref{progev} and similar arguments as above, we have
\begin{equation*}
\left|\cos\theta_\ell\right|_{G^{s_1,(1+\eta_{n-1})C}(I)},\ \  \left|\tan\theta_\ell\right|_{G^{s_1,(1+\eta_{n-1})C}(I)}\leq C,
\end{equation*}
\begin{equation*}
\left|\frac{1}{\cos\theta_\ell}\right|_{G^{s_1,(1+\eta_{n-1})C}(I)}, \ \ \left|\cot\theta\right|_{G^{s_1,(1+\eta_{n-1})C}(I)}\leq Ce^{q_n^{\gamma_1}}.
\end{equation*}
Thus all the assumptions in Lemma \ref{2} are satisfied, it follows
$$
\|A_{r_n^+}\|_{G^{s_1,(1+\eta_n)C}}\leq C^{2n_0}\prod_{\ell=0}^{n_0-1}|e^\ell|_{G^{s_1,(1+\eta_{n-1})K}}\leq C^{2n_0}\widetilde{\lambda_{n-1}}^{\sum_{\ell=0}^{n_0-1}(j_{\ell+1}-j_\ell)}\leq\widetilde{\lambda_n}^{r_n^+},
$$
$$
\left|\frac{1}{\|A_{r_n^+}\|}\right|_{G^{s_1,(1+\eta_i)C}}\leq C^{2n_0}e^{4n_0q_{n-1}^\gamma}\prod_{\ell=0}^{n_0-1}|(e^\ell)^{-1}|_{G^{s_1,(1+\eta_{n-1})K}}\leq \lambda_n^{-r_n^+},
$$
$$
|e_n|_{G^{s_1,(1+\eta_n)C}(I_n)}\leq \lambda_n^{-\frac{1}{20}q_{n-1}}.
$$
Similar results hold for $r_n^-$.

By the definition, we have
\begin{align*}
|\phi_n-\phi_{n-1}|_{G^{s_1,(1+\eta_n)C}(\mathbb{S}^1)}&=|\hat{e}_n|_{G^{s_1,(1+\eta_n)C}(\mathbb{S}^1)}\leq |e_n|_{G^{s_1,(1+\eta_n)C}(I_n)}|f_n|_{G^{s_1,(1+\eta_n)C}(\mathbb{S}^1)}\\
&\leq \lambda_{n-1}^{-\frac{q_{n-1}}{2}}(C_0q^\beta_n)^{q_n^{\frac{\nu_1\beta}{1-\delta_1\nu_1}}}\leq \lambda_n^{-\frac{1}{40}q_{n-1}}.
\end{align*}
The last inequality holds because  $\alpha$ is bounded.\\
\\

\noindent{\it Verifying $(2)_n$ of Proposition \ref{p1}}:  Define $A_n(x)=\Lambda\cdot R_{\frac{\pi}{2}-\phi_n(x)}$. Obviously, $A_n(x)=A_{n-1}(x)\cdot R_{-\hat{e}_n(x)}$.
\begin{Lemma}[\cite{wangyou}]
For $x\in I_n$, it holds that
$$
A_n^{r_n^+}(x)=A_{n-1}^{r_n^+}(x)\cdot R_{-\hat{e}_n(x)}
$$
and
$$
A_n^{-r_n^-}(x)=R_{\hat{e}_n(T^{-r_n^-}x)}\cdot A_{n-1}^{-r_n^-}(x).
$$
\end{Lemma}
Thus, for any $x\in I_n$, $\{A_n(x),...,A_n(T^{r_n^+(x)-1}x)\}$ is a $\lambda_n$-hyperbolic sequence.\\
\\

\noindent
{\it Verifying $(3)_n$ of Proposition \ref{p1}}: $(s_n-s_n^\prime)(x)=(\overline{s}_n-\overline{s}_n^\prime)(x)+\hat{e}_n(x)$ which implies that $(s_n-s_n^\prime)(x)=\phi_0(x)$ on $\frac{I_n}{10}$, since $|e_n(x)|_{G^{s,(1+\eta_n)C}}\leq \lambda_{n-1}^{-q_{n-1}}$ in $I_n$. Thus we have
$$
|(s_n-s_n^\prime)(x)|\geq |\phi_0(x)|-\lambda_{n-1}^{-q_{n-1}}\geq \frac{1}{2}e^{-(10q_n^\beta)^{\nu_1}},
$$
on $I_n\backslash\frac{I_n}{10}$.

Thus we finish the proof by letting $K_1=\lim\limits_{n\rightarrow\infty}(1+\eta_n)C$.
\subsection{Proof of Proposition \ref{p2}}
For any $n\geq N$, let $\tilde{e}_n(x)=-(s_n(x)-s'_n(x))\cdot f_n(x)$ be a $2\pi$-periodic smooth function such that it is $-(s_n(x)-s'_n(x))$ on $\frac{I_n}{10}$ and vanishes outside $I_n$. From $(3)_n$ in Proposition \ref{p1}, we have that $\tilde{e}_n(x)=-\phi_0(x)\cdot f_n(x)$. Then we define $\tilde{\phi}_n(x)=\phi_n(x)+\tilde{e}_n(x)$ and $\widetilde{A}_n(x)=\Lambda\cdot R_{\frac{\pi}{2}-\tilde{\phi}_n(x)}$.\\

\noindent {\it Verifying of $(1)_n$ of Proposition \ref{p2}}: It follows from the following Lemma.
\begin{Lemma}
$|\tilde{e}_n|_{G^{s_2,C}}\leq e^{-\frac{1}{4}q_n^{\gamma}}$ for $n>N$.
\end{Lemma}
\begin{pf}By \eqref{ind} in Lemma \ref{Gevrey-function}, we have $\|\phi_0\|_{G^{s_1,C}(I_n)}\leq e^{-\frac{1}{2 }q_n^{\beta\nu_1}}$ for some $C>0$. On the other hand,  we choose $\delta_2$ sufficiently small such that $\frac{\beta}{s_2-1-\delta_2}=\frac{\beta}{\frac{1}{\nu_2}-\delta_2}<\beta\nu_1$, by Lemma \ref{fn}, we have $$
|f_n|_{s_2,C}\leq (Cq^\beta_n)^{q_n^{\frac{\beta}{s_2-1-\delta_2}}}.
$$
Thus
$$
|\tilde{e}_n|_{G^{s_2,C}(\mathbb{S}^1)}\leq |\phi_0|_{G^{s_1,C}(I_n)}|f_n|_{G^{s_2,C}(\mathbb{S}^1)}\leq  e^{-\frac{1}{2 }q_n^{\beta\nu_1}}(Cq^\beta_n)^{q_n^{\frac{\beta}{s_2-1-\delta_2}}}\leq e^{-\frac{1}{4}q^{\gamma_1}_n}.
$$

\end{pf}
\noindent {\it Verifying  $(2)_n$ of Proposition \ref{p2}}: Since for each $x\in I_n$, $\{A_n(x), A_n(Tx),\cdots, A_n(T^{r_n^+(x)-1}x)\}$ is $\lambda_n$-hyperbolic and $\tilde{\phi}_n(x)=\phi_n(x)$ on $\mathbb{S}^1\backslash I_n$, we see that $\{\widetilde{A}_n(x), \widetilde{A}_n(Tx),\cdots, \widetilde{A}_n(T^{r_n^+(x)-1}x)\}$ is $\lambda_n$-hyperbolic. Thus $\tilde{s}_n(x)=s(\widetilde{A}^{r_n^+}_n(x))$ and $\tilde{s}'_n(x)=s(\widetilde{A}^{-r_n^+}_n(x))$ are well defined. \\
\\
\noindent
{\it Verifying $(3)_n$ of Proposition \ref{p2}}: Notice that $\tilde{s}_n(x)-\tilde{s}_n'(x)=s_n(x)-s_n'(x)-\tilde{e}_n(x)$.  Thus from the definition of $\tilde{e}_n(x)$, it holds that
\begin{align*}
\tilde{s}_n(x)=\tilde{s}_n^\prime(x)\ \  x\in \frac{I_n}{10}.
\end{align*}
Thus we finish the whole proof by choosing $\nu=\frac{1}{s_2-1}$ and $K=\max\{K_1,C\}$.

\section{The proofs of technical lemmas}
\subsection{Proof of Proposition \ref{progev}}The following two Lemmas will be used frequently.
\begin{Lemma}[The Formula of Faa di Bruno, see Theorem 1.3.2 in \cite{kp}]\label{1} Assume $f$ and $g$ are two smooth functions in an open interval $(a,b)$, let $h=g\circ f$, then
$$
h^{(n)}(x)=\sum_{k_1+2k_2+...+nk_n=n}\frac{n!}{k_1!k_2!...k_n!}g^{(k)}\left(f(x)\right)\left(\frac{f^{(1)}}{1!}\right)^{k_1}\left(\frac{f^{(2)}}{2!}\right)^{k_2}...\left(\frac{f^{(n)}}{n!}\right)^{k_n},
$$
where $k=k_1+k_2+...+k_n$.
\end{Lemma}
\begin{Lemma}[Lemma 1.4.1 in \cite{kp}]\label{com}
$$
\sum_{k_1+2k_2+...+nk_n=n}\frac{k!}{k_1!k_2!...k_n!}R^k=R(1+R)^{n-1},
$$
where $k=k_1+k_2+...+k_n$.
\end{Lemma}
By Stirling formula, one has
\begin{equation}\label{stirling}
\left(\frac{n}{e}\right)^n\leq n!\leq C\left(\frac{n}{e}\right)^n\sqrt{n}.
\end{equation}
\begin{Lemma}\label{use}
For any $0<\e\leq \frac{1}{2}$ and any $\sigma>0$, we have
$$
n^\sigma\leq (2\sigma)^\sigma\e^{-\sigma}(1+\e)^n.
$$
\end{Lemma}
\begin{pf}
Note that
$\frac{x^\sigma}{(1+\e)^x}=e^{-x\ln(1+\e)+\sigma\ln x}$. Let $f(x)=-x\ln(1+\e)+\sigma\ln x$, then
$$
f'(x)=-\ln(1+\e)+\frac{\sigma}{x}.
$$
It follows that $\max|f(x)|=f(\frac{\sigma}{\ln(1+\e)})=-\sigma+\sigma\ln\frac{\sigma}{\ln(1+\e)}$. Hence
$$
\frac{x^\sigma}{(1+\e)^x}\leq \left(\frac{\sigma}{\ln(1+\e)}\right)^\sigma\leq \left(\frac{2\sigma}{\e}\right)^\sigma.
$$
We finish the proof.
\end{pf}
\noindent {\bf Proof of Proposition \ref{progev}:} The proof of (1) and (2) can be found in \cite{bf}. Now we prove (3), note that $\left(\frac{1}{x}\right)^{(n)}=\frac{(-1)^nn!}{x^{n+1}}$. By Lemma \ref{1}, we have
\begin{equation}\label{df}
\left(\frac{1}{f}\right)^{(n)}=\sum_{k_1+2k_2+...+nk_n=n}\frac{n!}{k_1!k_2!...k_n!}\frac{(-1)^kk!}{f^{k+1}}\left(\frac{f^{(1)}}{1!}\right)^{k_1}\left(\frac{f^{(2)}}{2!}\right)^{k_2}...\left(\frac{f^{(n)}}{n!}\right)^{k_n},
\end{equation}
where $k=k_1+k_2+...+k_n$. Recall that
$$
|f|_{s,K}:=\frac{4\pi^2}{3}\sup\limits_{n}\frac{(1+|n|)^2}{K^n(n!)^s}|\partial^n f|_{C^0(I)},
$$
it follows that
\begin{equation}\label{df1}
\inf_{x\in I}\left|f(x)\right|\geq \frac{3}{4\pi^2}(1-\e),\ \ \sup_{x\in I}\left|f^{(n)}(x)\right|\leq \begin{cases} \frac{3}{4\pi^2}(1+\e) &n=0\\
\frac{3\e}{4\pi^2}\frac{K^n(n!)^{s}}{(1+n)^2}& n\geq 1
\end{cases}.
\end{equation}
Let $c= \frac{3}{4\pi^2}(1-\e)$. By \eqref{df} and  \eqref{df1}, for $n\geq 1$, we have
\begin{align*}
\left|\left(\frac{1}{f}\right)^{(n)}\right|&\leq c^{-1}\left|\sum_{k_1+2k_2+...+nk_n=n}\frac{n!}{k_1!k_2!...k_n!}\frac{k!}{c^{k}}\left(\frac{f^{(1)}}{1!}\right)^{k_1}\left(\frac{f^{(2)}}{2!}\right)^{k_2}...\left(\frac{f^{(n)}}{n!}\right)^{k_n}\right|\\
&\leq c^{-1}\left|\sum_{k_1+2k_2+...+nk_n=n}\frac{n!}{k_1!k_2!...k_n!}\frac{k!}{c^{k}}\e^k K^{n}\left(\frac{(2!)^{s}}{2!}\right)^{k_2}...\left(\frac{(n!)^{s}}{n!}\right)^{k_n}\right|.
\end{align*}
By \eqref{stirling} and Lemma \ref{use}, we have
$$
n!\leq C\left(\frac{n}{e}\right)^n\sqrt{n}\leq C\left(\frac{n}{e}\right)^n(1+\e^{\frac{1}{s+1}})^n\e^{-\frac{1}{2(s+1)}}.
$$
Thus
$$
\left|\left(\frac{(2!)^{s}}{2!}\right)^{k_2}...\left(\frac{(n!)^{s}}{n!}\right)^{k_n}\right|\leq C^{(s-1)k}\e^{-\frac{1}{2}k}\left(\frac{n}{e}\right)^{(s-1)n}(1+\e^{\frac{1}{s+1}})^{(s-1)n}.
$$
It follows that
\begin{align*}
\left|\left(\frac{1}{f}\right)^{(n)}\right|&\leq c^{-1}n! \left(\frac{n}{e}\right)^{(s-1)n}(1+\e^{\frac{1}{s+1}})^{(s-1)n}K^n\left|\sum_{k_1+2k_2+...+nk_n=n}\frac{k!}{k_1!k_2!...k_n!}(\e^{\frac{1}{2}} c^{-1}C^s)^{k}\right|\\
&\leq  \e^{\frac{1}{2}} c^{-2}C^s(n!)^s \left[K(1+\e^{\frac{1}{2}} c^{-1}C^s)(1+\e^{\frac{1}{s+1}})^{(s-1)}\right]^n,
\end{align*}
where  $k=k_1+k_2+...+k_n$ and the last inequality follows from Lemma \ref{com} and \eqref{stirling}.
For  sufficiently small $\e$ depending on $s$, we have
$$
(1+\e^{\frac{1}{2}} c^{-1}C^s)(1+\e^{\frac{1}{s-1}})^{(s-1)}\leq 1+\frac{1}{4}\e^{\frac{1}{s}},\ \ \e^{\frac{1}{2}} c^{-2}C^s\leq \e^{\frac{1}{3}}.
$$
Hence by Lemma \ref{use} again,
$$
\frac{4\pi^2}{3}\sup\limits_{n}\frac{(1+|n|)^2}{\left(K(1+\e^{\frac{1}{s+8}})\right)^n(n!)^s}\left|\partial^n \frac{1}{f}\right|_{C^0(I)}\leq \e^{\frac{1}{3}}|n|^2\left(1+\frac{1}{2}\e^{\frac{1}{s+8}}\right)^{-n}\leq \e^{\frac{1}{12}}.
$$
By the definition, we have
$$
\left|\frac{1}{f}-1\right|_{s,(1+\e^{\frac{1}{s+8}})K}\leq\e^{\frac{1}{12}}.
$$

For (4), note that $|(\sqrt{x})^n|=|\frac{1}{2}\cdots(\frac{1}{2}-n+1)x^{\frac{1}{2}-n}|\leq (n+2)!\sqrt{|x|}|x|^{-n}$. Similar to the proof of (2), we have
\begin{align*}
\left|\left(\sqrt{f}\right)^{(n)}\right|\leq  \e^{\frac{1}{2}} c^{-1}C^s(n+2)^2(n!)^s \left[K(1+\e^{\frac{1}{2}} c^{-1}C^s)(1+\e^{\frac{1}{s-1}})^{(s+1)}\right]^n.
\end{align*}
By Lemma \ref{use}, we have
$$
\frac{4\pi^2}{3}\sup\limits_{n}\frac{(1+|n|)^2}{\left(K(1+\e^{\frac{1}{s+16}})\right)^n(n!)^s}\left|\partial^n \sqrt{f}\right|_{C^0(I)}\leq \e^{\frac{1}{3}}|n|^4\left(1+\frac{1}{2}\e^{\frac{1}{s+16}}\right)^{-n}\leq \e^{\frac{1}{12}}.
$$

The proof of (5) is exactly the same as (2) since $|\arcsin^{(n)}x|\leq 2^nn!$, $|\sin^{(n)}x|,|\cos^{(n)}x|\leq 1\leq n!$ for any $n\in\N$. Thus we finish the proof.

\subsection{Proof of Lemma \ref{Gevrey-function}}
We inductively prove for $x>0$,
\begin{equation}\label{indss}
f^{(n)}(x)=\sum\limits_{i=1}^n\frac{a^n_i(\nu)}{x^{i\nu+n}}e^{-\frac{1}{x^\nu}},
\end{equation}
\begin{equation}\label{ind1}
|a^n_i(\nu)|\leq (2\nu+2)^{n+i}(\nu+n)^{n-i},\ \ 1\leq i\leq n.
\end{equation}
Assume for $k\leq n$, \eqref{indss} and \eqref{ind1} hold, then for $k=n+1$, we have
\begin{align*}
f^{(n+1)}(x)=\sum\limits_{i=1}^n\frac{\nu a^n_i(\nu)}{x^{i\nu+n+\nu+1}}e^{-\frac{1}{x^\nu}}-\sum\limits_{i=1}^n\frac{a^n_i(\nu)(i\nu+n)}{x^{i\nu+n+1}}e^{-\frac{1}{x^\nu}}:=\sum\limits_{i=1}^{n+1}\frac{a^{n+1}_i(\nu)}{x^{i\nu+n+1}}e^{-\frac{1}{x^\nu}},
\end{align*}
where
$$
a^{n+1}_{i}=\left\{
\begin{aligned}
&-a_1^n(\nu)(\nu+n)& i=1\\
&a_{i-1}^n(\nu)\nu-a_{i}^n(\nu)(i\nu+n)& 2\leq i\leq n\\
&a_n^n(\nu)\nu& i=n+1\\
\end{aligned}
\right. .
$$
By \eqref{ind1}, we have
$$
|a^{n+1}_{i}|\leq \left\{
\begin{aligned}
&|a_i^n(\nu)|(\nu+n)\leq (2\nu+2)^{n+i}(\nu+n)^{n+1-i}\leq (2\nu+2)^{n+1+i}(\nu+n+1)^{n+1-i} & i=1\\
&|a_{i-1}^n(\nu)\nu|+|a_{i}^n(\nu)(i\nu+n)|\leq(2\nu+2)^{n+i+1}(\nu+n+1)^{n+1-i} & 2\leq i\leq n\\
&|a_i^n(\nu)\nu|\leq  (2\nu+2)^{n+1+i}(\nu+n)^{n-i}\leq (2\nu+2)^{n+1+i}(\nu+n+1)^{n+1-i}& i=n+1\\
\end{aligned}
\right. .
$$
\eqref{indss} and \eqref{ind1} imply that
\begin{equation}\label{f_deri}
|f^{(n)}(x)|\leq \sum\limits_{i=1}^n\frac{(2\nu+2)^{n+i}(\nu+n)^{n-i}}{|x|^{i\nu+n}}e^{-\frac{1}{|x|^\nu}}.
\end{equation}
Notice that
$$
\sup\limits_{x\in\R}\frac{1}{|x|^{i\nu+n}}e^{-\frac{1}{2|x|^\nu}}\leq \left(\frac{2(i\nu+n)}{\nu}\right)^{i+n/\nu},
$$
it follows that
\begin{align}\label{ind}
\nonumber |f^{(n)}(x)|&\leq e^{-\frac{1}{2|x|^\nu}}\sum\limits_{i=1}^n(2\nu+2)^{n+i}(\nu+n)^{n-i}\left(\frac{2(i\nu+n)}{\nu}\right)^{i+n/\nu}\\
&\leq e^{-\frac{1}{2|x|^\nu}}C^{n}(\nu+n)^{n(1+\frac{1}{\nu})}\leq e^{-\frac{1}{2|x|^\nu}}C^n(n!)^{1+\frac{1}{\nu}}.
\end{align}
\subsection{Proof of Corollary \ref{inverse-Gevrey-function}}By the same argument as in Lemma \ref{Gevrey-function}, we have for any $x> 0$,
\begin{equation}\label{ind'}
f^{(n)}(x)=\sum\limits_{i=1}^n\frac{a^n_i(\nu)}{x^{i\nu+n}}e^{\frac{1}{x^\nu}},
\end{equation}
\begin{equation}\label{ind1'}
|a^n_i(\nu)|\leq (2\nu+2)^{n+i}(\nu+n)^{n-i},\ \ 1\leq i\leq n.
\end{equation}
\eqref{ind'} and \eqref{ind1'} imply that
\begin{equation}\label{f_deri'}
|f^{(n)}(x)|\leq \sum\limits_{i=1}^n\frac{(2\nu+2)^{n+i}(\nu+n)^{n-i}}{|x|^{i\nu+n}}e^{\frac{1}{|x|^\nu}}.
\end{equation}
Notice that
$$
\sup\limits_{x\in\R}\frac{1}{|x|^{i\nu+n}}e^{-\frac{1}{|x|^\nu}}\leq \left(\frac{i\nu+n}{\nu}\right)^{i+n/\nu},
$$
it follows that
$$
|f^{(n)}(x)|\leq e^{\frac{2}{|x|^\nu}} \sum\limits_{i=1}^n(2\nu+2)^{n+i}(\nu+n)^{n-i}\left(\frac{i\nu+n}{\nu}\right)^{i+n/\nu}\leq e^{\frac{2}{|x|^\nu}}C^n(n!)^{1+\frac{1}{\nu}}.
$$
\subsection{Proof of Corollary \ref{phi_0}}
For any $k\in\Z$ and $x\in [c_1+k\pi,c_1+(k+1)\pi)$, we have
$$
g(x)=cf(x-c_1-k\pi)f(c_1+(k+1)\pi-x).
$$
Let us firstly show that $g\in C^\infty(\mathbb{S}^1)$, for which we only need to verify the derivative exists for $x=c_1+k\pi$. By a direct calculation, we have
$$
g^{(n)}(x)=\begin{cases}\sum\limits_{k=0}^n\begin{pmatrix}n\\ k\end{pmatrix} f^{(n-k)}(x-c_1-k\pi)(-1)^kf^{(k)}(c_1+(k+1)\pi-x)&c_1+k\pi<x<c_1+(k+1)\pi\\
\sum\limits_{k=0}^n\begin{pmatrix}n\\ k\end{pmatrix} f^{(n-k)}(x-c_1-(k-1)\pi)(-1)^kf^{(k)}(c_1+k\pi-x)& c_1+(k-1)\pi<x<c_1+k\pi
\end{cases}.
$$
We inductively prove that
\begin{equation}\label{deri}
g_+^{(n)}(c_1+k\pi)=g_-^{(n)}(c_1+k\pi)=0.
\end{equation}
Assume \eqref{deri} holds for all $k\leq n$. For $k=n+1$, by Lemma \ref{Gevrey-function}, we have
\begin{align*}
&\lim\limits_{x\searrow c_1+k\pi}\left|\frac{g^{(n)}(x)-g^{(n)}(c_1+k\pi)}{x-c_1-k\pi}\right|\\
\leq  &\lim\limits_{x\searrow c_1+k\pi}\frac{\sum\limits_{k=0}^n\begin{pmatrix}n\\ k\end{pmatrix}e^{-\frac{1}{2|x-c_1-k\pi|^\nu}}C^{n-k}((n-k)!)^{1+\frac{1}{\nu}}e^{-\frac{1}{2|c_1+(k+1)\pi-x|^\nu}}C^{k}((k)!)^{1+\frac{1}{\nu}}}{x-c_1-k\pi}=0.
\end{align*}
\begin{align*}
&\lim\limits_{x\nearrow c_1+k\pi}\left|\frac{g^{(n)}(x)-g^{(n)}(c_1+k\pi)}{x-c_1-k\pi}\right|\\
\leq  &\lim\limits_{x\nearrow c_1+k\pi}\frac{\sum\limits_{k=0}^n\begin{pmatrix}n\\ k\end{pmatrix}e^{-\frac{1}{2|x-c_1-(k-1)\pi|^\nu}}C^{n-k}((n-k)!)^{1+\frac{1}{\nu}}e^{-\frac{1}{2|c_1+k\pi-x|^\nu}}C^{k}((k)!)^{1+\frac{1}{\nu}}}{|x-c_1-k\pi|}=0.
\end{align*}
Thus
\begin{equation*}
g_+^{(n+1)}(c_1+k\pi)=g_-^{(n+1)}(c_1+k\pi)=0.
\end{equation*}
By (1) in Proposition \ref{progev} and Lemma \ref{Gevrey-function}, we have
$$
|g|_{s,C}\leq |f(\cdot-c_1)|_{s,C}|f(c_1+\pi-\cdot)|_{s,C}\leq ce^{-\frac{1}{2|x-c_1|^\nu}}e^{-\frac{1}{2|x-c_1-\pi|^\nu}}.
$$
We thus finish the proof.

\subsection{Proof of Lemma \ref{fn}}
Define
$$
\phi(x)=
\begin{cases}
e^{-\frac{1}{x^{1/\delta}}}&x>0\\
0&x\leq 0
\end{cases}.
$$
Let
$$
w_1(x)=
\begin{cases}
w_0(-x)&x>0\\
w_0(x)&x\leq 0
\end{cases},
$$
where $w_0(x)=\frac{\phi(x+2)}{\phi(x+2)+\phi(-x-1)}$.  It's easy to verify that
$$
w_1(x)=\left\{
\begin{aligned}
&0& x\geq 2\\
&\frac{e^{-\frac{1}{(-x+2)^{1/\delta}}}}{e^{-\frac{1}{(-x+2)^{1/\delta}}}+e^{-\frac{1}{(x-1)^{1/\delta}}}}& 1< x<2\\
&1& -1\leq x\leq 1\\
&\frac{e^{-\frac{1}{(x+2)^{1/\delta}}}}{e^{-\frac{1}{(x+2)^{1/\delta}}}+e^{-\frac{1}{(-x-1)^{1/\delta}}}}& -2< x<-1\\
&0& x\leq -2
\end{aligned}
\right. .
$$
By similar arguments as Corollary \ref{phi_0}, we have $w_1\in G^{1+{\delta}}(\R)$.

Then we define $f_n$ to be a $\pi$-periodic function such that
$$
f_n(x)=w_1(10q_n^\beta(x-c_1)), \ \ x\in \left[c_1-\frac{\pi}{2},c_1+\frac{\pi}{2}\right].
$$
From the definition, we have that $f_n^{(r)}(x)=(10q_n)^{\beta r}\cdot w^{(r)}_1(y)$ where $y=10q_n^2(x-c_1)$.  By the definition of $G^{1+{\delta}}$-norm, there exists $C>0$ such that
$$
\sup\limits_{x\in I_n}|f_n^{(r)}(x)|\leq \frac{(Cq_n)^{\beta r}(r!)^{1+\delta}}{1+r^2}.
$$
Thus
$$
\sup\limits_{x\in I_n}\frac{|f_n^{(r)}(x)|(1+r^2)}{C^r(r!)^{1+\frac{1}{\nu}}}\leq q_n^{\beta r}(r!)^{\delta-\frac{1}{\nu}}\leq (Cq^\beta_n)^{q_n^{\frac{\nu\beta}{1-\delta\nu}}}.
$$
\section*{Acknowledgement}
We would like to thank Svetlana Jitomirskaya for many valuable discussions. L. Ge and X. Zhao were partially supported by NSF DMS-1901462. L. Ge was partially supported by AMS-Simons Travel Grant 2020-2022. Y. Wang was supported by NSF of China (11771205). J. You was partially supported by National Key R\&D Program of China (2020YFA0713300), NNSF of China (11871286).  X. Zhao was partially supported by China Scholarship Council (No. 201906190072).

\end{document}